\documentclass[leqno,11pt]{amsart}

\usepackage{graphicx}
\usepackage{latexsym}
\usepackage{amssymb}
\usepackage{psfrag}
\usepackage{eucal}
\usepackage{pb-diagram}
\usepackage[bottom=1.0in, right=1.0in, top=1.6in, left=1.4in, marginparsep=0in, marginparwidth=0in, headheight=0in, headsep=0.3in, footskip=0.3in, includefoot]{geometry}
\newcommand\nc{\newcommand}
\nc\rnc{\renewcommand}
\nc\oo{^{-1}}
\nc\half{{\frac12}}
\nc\Mod{{\mathrm{Mod}}}
\nc\tr{{\mathrm{tr}}}
\nc\Id{{\mathrm{Id}}}
\nc\ad{{\mathrm{ad}}}
\nc\inn{{\mathrm{in}}}
\nc\ifi{{\mathrm{if}}}
\nc\Hom{{\mathrm{Hom}}}
\nc\Ker{{\mathrm{Ker}}}
\nc\End{{\mathrm{End}}}
\nc\Aut{{\mathrm{Aut}}}
\nc\cen{{\mathrm{center}}}
\nc\sn{{\mathrm{sn}}}
\nc\dij{{\delta_{i j}}}
\nc\st{{\ |\ }}
\nc\FMod{{\mathrm{FMod}}}
\nc\Flip{{\mathrm{Flip}}}
\nc\im{{i = 1, \ldots, m}}
\nc\jn{{j = 1, \ldots, n}}
\nc\obj{{\mathrm{obj}}}
\def\be{\begin{equation}}
\def\ee{\end{equation}}
\def\ba{\begin{array}}
\def\ea{\end{array}}
\def\bea{\begin{eqnarray}}
\def\eea{\end{eqnarray}}
\nc\cA{{\mathcal A}}
\nc\cB{{\mathcal B}}
\nc\cC{{\mathcal C}}
\nc\cD{{\mathcal D}}
\nc\cE{{\mathcal E}}
\nc\cF{{\mathcal F}}
\nc\cG{{\mathcal G}}
\nc\cH{{\mathcal H}}
\nc\cI{{\mathcal I}}
\nc\cJ{{\mathcal J}}
\nc\cK{{\mathcal K}}
\nc\cL{{\mathcal L}}
\nc\cM{{\mathcal M}}
\nc\cN{{\mathcal N}}
\nc\cO{{\mathcal O}}
\nc\cP{{\mathcal P}}
\nc\cQ{{\mathcal Q}}
\nc\cR{{\mathcal R}}
\nc\cS{{\mathcal S}}
\nc\cT{{\mathcal T}}
\nc\cU{{\mathcal U}}
\nc\cV{{\mathcal V}}
\nc\cW{{\mathcal W}}
\nc\cX{{\mathcal X}}
\nc\cY{{\mathcal Y}}
\nc\cZ{{\mathcal Z}}
\nc\bA{{\mathbf{A}}}
\nc\bB{{\mathbf{B}}}
\nc\bC{{\mathbf{C}}}
\nc\bD{{\mathbf{D}}}
\nc\bE{{\mathbf{E}}}
\nc\bF{{\mathbf{F}}}
\nc\bG{{\mathbf{G}}}
\nc\bH{{\mathbf{H}}}
\nc\bI{{\mathbf{I}}}
\nc\bJ{{\mathbf{J}}}
\nc\bK{{\mathbf{K}}}
\nc\bL{{\mathbf{L}}}
\nc\bM{{\mathbf{M}}}
\nc\bN{{\mathbf{N}}}
\nc\bO{{\mathbf{O}}}
\nc\bP{{\mathbf{P}}}
\nc\bQ{{\mathbf{Q}}}
\nc\bR{{\mathbf{R}}}
\nc\bS{{\mathbf{S}}}
\nc\bT{{\mathbf{T}}}
\nc\bU{{\mathbf{U}}}
\nc\bV{{\mathbf{V}}}
\nc\bW{{\mathbf{W}}}
\nc\bX{{\mathbf{X}}}
\nc\bY{{\mathbf{Y}}}
\nc\bZ{{\mathbf{Z}}}
\nc\BA{{\mathbb{A}}}
\nc\BB{{\mathbb{B}}}
\nc\BC{{\mathbb{C}}}
\nc\BD{{\mathbb{D}}}
\nc\BE{{\mathbb{E}}}
\nc\BF{{\mathbb{F}}}
\nc\BG{{\mathbb{G}}}
\nc\BH{{\mathbb{H}}}
\nc\BI{{\mathbb{I}}}
\nc\BJ{{\mathbb{J}}}
\nc\BK{{\mathbb{K}}}
\nc\BL{{\mathbb{L}}}
\nc\BM{{\mathbb{M}}}
\nc\BN{{\mathbb{N}}}
\nc\BO{{\mathbb{O}}}
\nc\BP{{\mathbb{P}}}
\nc\BQ{{\mathbb{Q}}}
\nc\BR{{\mathbb{R}}}
\nc\BS{{\mathbb{S}}}
\nc\BT{{\mathbb{T}}}
\nc\BU{{\mathbb{U}}}
\nc\BV{{\mathbb{V}}}
\nc\BW{{\mathbb{W}}}
\nc\BX{{\mathbb{X}}}
\nc\BY{{\mathbb{Y}}}
\nc\BZ{{\mathbb{Z}}}
\nc\Bone{{\mathbb{1}}}
\nc\bb{{\mathbf{b}}}
\nc\bc{{\mathbf{c}}}
\nc\bd{{\mathbf{d}}}
\nc\bbf{{\mathbf{f}}}
\nc\bg{{\mathbf{g}}}
\nc\bh{{\mathbf{h}}}
\nc\bi{{\mathbf{i}}}
\nc\bj{{\mathbf{j}}}
\nc\bk{{\mathbf{k}}}
\nc\bl{{\mathbf{l}}}
\nc\bm{{\mathbf{m}}}
\nc\bn{{\mathbf{n}}}
\nc\bo{{\mathbf{o}}}
\nc\bp{{\mathbf{p}}}
\nc\bq{{\mathbf{q}}}
\nc\br{{\mathbf{r}}}
\nc\bs{{\mathbf{s}}}
\nc\bt{{\mathbf{t}}}
\nc\bu{{\mathbf{u}}}
\nc\bv{{\mathbf{v}}}
\nc\bw{{\mathbf{w}}}
\nc\bx{{\mathbf{x}}}
\nc\by{{\mathbf{y}}}
\nc\bz{{\mathbf{z}}}
\nc\0{\mathbf{0}}
\nc\1{\mathbf{1}}
\nc\fa{{\mathfrak{a}}}
\nc\fb{{\mathfrak{b}}}
\nc\fc{{\mathfrak{c}}}
\nc\fd{{\mathfrak{d}}}
\nc\fe{{\mathfrak{e}}}
\nc\ff{{\mathfrak{f}}}
\nc\fg{{\mathfrak{g}}}
\nc\fh{{\mathfrak{h}}}
\nc\fj{{\mathfrak{j}}}
\nc\fk{{\mathfrak{k}}}
\nc\fl{{\mathfrak{l}}}
\nc\fm{{\mathfrak{m}}}
\nc\fn{{\mathfrak{n}}}
\nc\fo{{\mathfrak{o}}}
\nc\fp{{\mathfrak{p}}}
\nc\fq{{\mathfrak{q}}}
\nc\fr{{\mathfrak{r}}}
\nc\fs{{\mathfrak{s}}}
\nc\ft{{\mathfrak{t}}}
\nc\fu{{\mathfrak{u}}}
\nc\fv{{\mathfrak{v}}}
\nc\fw{{\mathfrak{w}}}
\nc\fx{{\mathfrak{x}}}
\nc\fy{{\mathfrak{y}}}
\nc\fz{{\mathfrak{z}}}
\nc\fA{{\mathfrak{A}}}
\nc\fB{{\mathfrak{B}}}
\nc\fC{{\mathfrak{C}}}
\nc\fD{{\mathfrak{D}}}
\nc\fE{{\mathfrak{E}}}
\nc\fF{{\mathfrak{F}}}
\nc\fG{{\mathfrak{G}}}
\nc\fH{{\mathfrak{H}}}
\nc\fI{{\mathfrak{I}}}
\nc\fJ{{\mathfrak{J}}}
\nc\fK{{\mathfrak{K}}}
\nc\fL{{\mathfrak{L}}}
\nc\fM{{\mathfrak{M}}}
\nc\fN{{\mathfrak{N}}}
\nc\fO{{\mathfrak{O}}}
\nc\fP{{\mathfrak{P}}}
\nc\fQ{{\mathfrak{Q}}}
\nc\fR{{\mathfrak{R}}}
\nc\fS{{\mathfrak{S}}}
\nc\fT{{\mathfrak{T}}}
\nc\fU{{\mathfrak{U}}}
\nc\fV{{\mathfrak{V}}}
\nc\fW{{\mathfrak{W}}}
\nc\fX{{\mathfrak{X}}}
\nc\fY{{\mathfrak{Y}}}
\nc\fZ{{\mathfrak{Z}}}
\nc\ga{{\alpha}}
\nc\gb{{\beta}}
\nc\gd{{\delta}}
\nc\gge{{\epsilon}}
\nc\gve{{\varepsilon}}
\nc\gz{{\zeta}}
\nc\gh{{\eta}}
\nc\gy{{\theta}}
\nc\gvy{{\vartheta}}
\nc\gi{{\iota}}
\nc\gk{{\kappa}}
\nc\gl{{\lambda}}
\nc\gm{{\mu}}
\nc\gn{{\nu}}
\nc\gx{{\xi}}
\nc\go{{o}}
\nc\gp{{\pi}}
\nc\gvp{{\varpi}}
\nc\gr{{\rho}}
\nc\gs{{\sigma}}
\nc\gvs{{\varsigma}}
\nc\gt{{\tau}}
\nc\gu{{\upsilon}}
\nc\gf{{\phi}}
\nc\gvf{{\varphi}}
\nc\emp{{\gf}}
\nc\gq{{\psi}}
\nc\gw{{\omega}}
\nc\gc{{\chi}}
\nc\gG{{\Gamma}}
\nc\gD{{\Delta}}
\nc\gY{{\Theta}}
\nc\gL{{\Lambda}}
\nc\gX{{\Xi}}
\nc\gP{{\Pi}}
\nc\gU{{\Upsilon}}
\nc\gF{{\Phi}}
\nc\gQ{{\Psi}}
\nc\gW{{\Omega}}
\nc\gS{{\Sigma}}
\newtheorem{thm}{{Theorem}}[section]
\newtheorem{lemma}[thm]{{Lemma}}
\newtheorem{prop}[thm]{{Proposition}}
\newtheorem{cor}[thm]{{Corollary}}
\newtheorem{dfn}[thm]{{Definition}}

\newtheorem{eg}[thm]{{Example}}

\nc{\brmk}{\vspace{2ex}\noindent{\em Remark.\ }}
\nc{\ermk}{\vspace{2ex}}
\nc{\beg}{\begin{eg}\em}
\nc{\eeg}{$\hfill\Box$\end{eg}}
\nc{\bdfn}{\begin{dfn}\em}
\nc{\edfn}{\end{dfn}}
\nc{\bproof}{{

\noindent {\em Proof. }}}
\nc{\eproof}{{$\hfill\Box$

\vspace*{1em}

}}
\nc{\thmref}[1]{theorem~\ref{#1}}
\nc{\secref}[1]{section~\ref{#1}}
\nc{\lemmaref}[1]{lemma~\ref{#1}}
\nc{\figref}[1]{figure~\ref{#1}}
\nc{\corref}[1]{corollary~\ref{#1}}
\nc{\factref}[1]{fact~\ref{#1}}
\nc{\propref}[1]{proposition~\ref{#1}}
\nc{\dfnref}[1]{definition~\ref{#1}}
\nc{\conref}[1]{conjecture~\ref{#1}}
\nc{\rmkref}[1]{remark~\ref{#1}}
\nc{\tabref}[1]{table~\ref{#1}}
\nc{\egref}[1]{example~\ref{#1}}
\nc{\chapref}[1]{chapter~\ref{#1}}

\nc\rrp{\hbox{\rm )}}
\nc\rlp{\hbox{\rm (}}
\nc{\eqn}[1]{{\rlp\ref{#1}\rrp}}
\nc\intt{\hbox{\rm int}}
\nc\supp{\hbox{\rm supp}}
\nc\il{{i = 1, \ldots, \ell}}
\nc\qv{{\BQ(v)}}
\nc\qx{{\BQ(\xi)}}
\nc\zx{{\BZ[\xi]}}
\nc\zz{{\BZ[\zeta]}}
\nc\zxr{{\BZ[\xi,r\oo]}}
\nc\zzr{{\BZ[\zeta,r\oo]}}
\nc\zv{{\BZ[v, v\oo]}}
\nc\An{{\cA_n}}
\nc\Ainf{{\cA_\infty}}
\nc\UA{\mathbf{U}\!_\cA}
\nc\Uv{\mathbf{U}_v}
\nc\UC{\mathbf{U}_\BC}
\nc\Ux{\mathbf{U}_\xi}
\nc\gVl{{_\fg\! V_\lambda}}
\nc\vVl{{_v\! V_\lambda}}
\nc\vV{{_v\! V}}
\nc\AVl{{_\cA\! V_\lambda}}
\nc\xVl{{_\xi\! V_\lambda}}
\nc\xVmi{{_\xi\! V_{\mu_i}}}
\nc\xVli{{_\xi\! V_{\gl_i}}}
\nc\xVlo{{_\xi\! V_{\gl_1}}}
\nc\xVls{{_\xi\! V_{\gl_6}}}
\nc\xVlsv{{_\xi\! V_{\gl_7}}}
\nc\xVmo{{_\xi\! V_{\gm_1}}}
\nc\xVmm{{_\xi\! V_{\gm_m}}}
\nc\xVm{{_\xi\! V_\mu}}
\nc\zVl{{_\gz\! V_\lambda}}
\nc\zVli{{_\gz\! V_{\lambda_i}}}
\nc\zVm{{_\gz\! V_\mu}}
\nc\zVmo{{_\gz\! V_{\mu_1}}}
\nc\zVmm{{_\gz\! V_{\mu_m}}}
\nc\xVwl{{_\xi\! V_{\lambda'}}}
\nc\CVl{{_\BC\! V_\lambda}}
\nc\AnVl{{_\An\!\!\! V_\lambda}}
\nc\AVn{{_\cA\! V_n}}
\nc\AnVm{{_\An\!\!\! V_\mu}}
\nc\AnVwl{{_\An\!\!\! V_{\lambda'}}}
\nc\AinfVl{{_\Ainf\!\!\!\! V_\lambda}}
\nc\AinfV{_\Ainf\!\!\!\! V}
\nc\Vv{{\cV_v}}
\nc\VC{\cV_\BC}
\nc\Vx{\cV_\xi}
\nc\Vz{\cV_\zeta}
\nc\bVx{\bar{\cV}_\xi}
\nc\bVz{\bar{\cV}_\zeta}
\nc\Vg{{\cV_\fg}}
\nc\VA{{\cV_{\!\cA}}}
\nc\CR{{\cC(\cR)}}

\begin{document}

\title{On certain integral tensor categories and integral TQFTs}
\author{Qi Chen}
\email{qichen@buffalo.edu}

\begin{abstract}
We construct certain tensor categories that are dominated by finitely many simple objects. Objects
in these categories are modules over rings of algebra integers. We show how to 
obtain TQFTs defined over algebra integers
from these categories.
\end{abstract}

\maketitle

%%%%%%%%%%%%%%%%%%%%%%%%%%%%%%%%%%%%%%%%%%%%%%%%%%%%%%%%
%%%%%%%%%%%%%%%%%%%%%%%%%%%%%%%%%%%%%%%%%%%%%%%%%%%%%%%%
\section{Introduction}\label{intro}
The goal of this note is twofold: (1) to give certain tensor categories 
dominated by finitely many simple objects (cf. \secref{sec:tensorDef}), 
(2) to construct topological quantum field theories (TQFTs) whose ground rings
are algebraic integers (or integral TQFTs for short) (cf. \secref{tqfts}).
We use representation theory of quantum groups in the construction. Results
in this note are related to those in \cite{andersen}, \cite{gk} and \cite{g2}.

We fix a complex simple Lie algebra $\fg$, an integer $r$ and a primitive $r$-th root of unity $\xi$
throughout this note.

Let $\Uv = \Uv(\fg)$ be the quantum group associated to $\fg$. It is well known that
$\Uv$ is a Hopf algebra over the ring $\qv$ where $v$ is a formal parameter.
Hence the category $\Vv$ of finite-dimensional $\Uv$-modules of type 1 is a tensor category (cf. \secref{sec:quantum}). 
It turns out that $\Vv$ is equivalent to the category 
of finite-dimensional $\fg$-modules (cf. 6.3 of \cite{lu}).

Let $\cA = \zv$. Lusztig showed that $\Uv$ has an $\cA$-subalgebra 
$\UA$ which inherits the Hopf algebra structure from $\Uv$. Let
$\UC = \UA\otimes_\cA \BC$ (resp. $\Ux = \UA\otimes_\cA \zx$) where $\BC$ (resp. $\zx$)
is considered as an $\cA$-algebra by sending $v$ to $\xi$. The category 
$\VC$ (resp. $\Vx$) of finite-dimensional $\UC$-modules (resp. finite-ranked $\Ux$-modules) 
of type 1 with the ground
ring $\BC$ (resp. $\zx$) behaves very differently from $\Vv$. 

If $r$ is prime to the non-zero entries of the Cartan matrix for $\fg$ and is bigger than the Coxeter number,
then Andersen proved in \cite{andersen} that (in the language of \cite{kirillov}) 
there is a full subcategory $\VC'$ of $\VC$ which has a
quotient category dominated by finitely many simple objects. (This is proved in \cite{gk} 
when $r$ is prime and big enough.)
Objects in these categories are vector spaces over $\BC$. The following proposition is our first main result.

\begin{prop}\label{prop:tensor}
If $\fg$ is not of type $E_8, F_4$ or $G_2$ and $r$ is a prime bigger than $m(\fg)$
then there exists a subcategory $\Vx'$ of $\Vx$ which has a quotient category $\bVx'$
dominated by finitely many simple objects.
\end{prop}

See \secref{root1} for $m(\fg)$ and the construction for the categories.
Note that objects in $\bVx'$ are $\zx$-modules of finite rank. 
This is essential to our construction of integral 
TQFTs. In the construction we have to restrict the allowed 3-cobordisms (cf. \secref{int_tqft}). 
TQFTs for such restricted 
3-cobordisms will be denoted by $(\cS_+, \tau)$ (as opposed to $(\cT, \tau)$ for the non-restricted ones).

\begin{prop}\label{prop:TQFT}
If $r$ and $\fg$ are as in \propref{prop:tensor}
then there exists a TQFT $(\cS_+, \tau)$ such that
the images of $\cS_+$ are free $\zx$-modules of finite rank.
\end{prop}

This proposition is contained in \propref{prop:integral}.
It provides us infinitely many mapping class group representations over $\zx$, 
which may shed lights on the structure of mapping class groups.

The note is organized as follows. Basic knowledge is recalled in \secref{pre}. Sections~\ref{rep} and
\ref{root1} deal with representations of quantum groups at generic parameter and root of unity respectively.
Proof of \propref{prop:decomposition} occupies \secref{sec:proof}.
Integral TQFTs are constructed in \secref{tqft}.

%%%%%%%%%%%%%%%%%%%%%%%%%%%%%%%%%%%%%%
%%%%%%%%%%%%%%%%%%%%%%%%%%%%%%%%%%%%%%
\section{Preliminaries}\label{pre}
The material in this section is quite standard and may be skipped.

%%%%%%%%%%%%%%%%%%%%%%%%%%%%%%%%%%%%%%
\subsection{Tensor categories}\label{sec:tensorDef}
A good reference for this subsection is \cite{kassel}. A category $\cC$ is a tensor category if
\begin{enumerate}
	\item[(a)] each hom-set $\Hom_\cC(U, V)$ is an additive abelian group and composition of morphisms is bilinear
	relative to this addition,
	\item[(b)] $\cC$ has biproduct (or direct sum) $U \oplus V$ for any two objects $U$ and $V$,
	\item[(c)] there is a bifunctor $\otimes : \cC \times \cC \to \cC$ which is associative up to a natural isomorphism,
	called the associativity constraint,
	$$
	\alpha_{U,V,W}: U\otimes(V\otimes W) \to (U\otimes V)\otimes W,
	$$
	\item[(d)] $\cC$ has an object $\1_\cC$, unique up to natural isomorphism, and two natural isomorphisms,
	called the left and the right unit constraint,
	$$
	\lambda_U : \1_\cC\otimes U \to U, \quad \rho_U : U \otimes \1_\cC \to U,
	$$
	\item[(e)] associativity, left unit and right unit constraints satisfy the pentagon and the triangle axioms.
\end{enumerate}

Our definition of a tensor category is different from the usual one, which assumes neither (a) nor (b).
Let $\cC$ be a tensor category. We say $\cC$ is dominated by a set of objects $A$ if 
every object is a direct sum of objects from $A$. Such a set is called a dominant set of $\cC$.
The ground ring of $\cC$ is $k_\cC = \Hom_\cC(\1_\cC, \1_\cC)$.
An object $U$ in $\cC$ is simple if $\Hom_\cC(U,U) = k_\cC$.
If $\cC$ is dominated by the set of simple objects then it is said to be semisimple.

A tensor category is said to be strict if its associativity, left unit and right unit constraints are all identities.
By Mac Lane's coherence theorem, every tensor category is equivalent to a strict one. Therefore we assume all tensor
categories in this note to be strict. A tensor category $\cC$ is said to have duality if
\begin{enumerate}
	\item[(f)] there is a contravariant functor 
	$* : \cC \to \cC$, (For any object $U$ and any morphism $f$ we denote $*(U)$ and
	$*(f)$ by $U^*$ and $f^*$ respectively.)
	\item[(g)] for every object $U$ there are two morphisms, called the coevaluation and the evaluation respectively,
	$$
	b_U : \1_\cC \to U\otimes U^*,\quad d_U: U^*\otimes U\to \1_\cC,
	$$
	such that
	$$
	(\Id_U\otimes d_U)(b_U\otimes \Id_U) = \Id_U,\quad (d_U\otimes\Id_{U^*})(\Id_{U^*}\otimes b_U) = \Id_{U^*}.
	$$
\end{enumerate}

Let $\cC$ be a tensor category and $P: \cC \times \cC \to \cC\times \cC$ be the flip functor defined by $P(U,V)=(V,U)$.
Then $\cC$ is said to be braided if there is a natural isomorphism $c$ between $\otimes$ and $\otimes P$ which
satisfies the hexagon axiom. Denote $c_{U,V}:U\otimes V\to V\otimes U$ for the braiding between $U$ and $V$.

Let $\cC$ be a braided tensor category with duality. A twist is a natural transformation $\theta$ from the identity
functor of $\cC$ to itself such that $\theta_{V\otimes W} = (\theta_V\otimes \theta_W)c_{W,V}c_{V,W}$ and
$\theta_{V^*} = \theta_V^*$ for all objects $V, W$ in $\cC$. 
A braided tensor category with duality is called a ribbon category if it has a twist.

Let $\cC$ be a ribbon category with a finite dominant set $A$. Set
$s_{U,V} = \tr_q(c_{U,V}c_{V,U})$ for every pair $U, V\in A$. 
Here 
$$
\tr_q(f) = \tr(d_Wc_{W,W^*}(f\theta_W\otimes \Id)b_W) \in k_\cC
$$
is the quantum trace of $f :W\to W$. The matrix $(s_{U,V})_{U,V\in A}$ is called the $S$-matrix of $\cC$ with respect
to the dominant set $A$ and is denoted by $S(\cC, A)$ or simply $S(\cC)$ if $A$ is clear from the context.

A ribbon category $\cC$ is said to be a modular category if
\begin{enumerate}
	\item [(h)] it has a finite dominant set $A$,
	\item [(i)] $A$ contains the identity object $\1_\cC$ and is closed under duality,
	\item [(j)] $S(\cC, A)$ is invertible over $k_\cC$.
\end{enumerate}

%%%%%%%%%%%%%%%%%%%%%%%%%%%%%%%%%%%%%%%%%%%%%%%%
\subsection{Lie algebras}\label{sec:Lie}
Let ($a_{ij}$), $i, j = 1, \ldots, \ell$, be the Cartan matrix and $\Phi$ be the root system of $\fg$. 
Fix a Cartan subalgebra $\fh$ of $\fg$ and a set of simple roots 
$\Pi_\fh = \{\alpha_1, \ldots, \alpha_\ell\}$ in its dual space $\fh^*$. We consider $\Pi_\fh$ as a subset of $\Phi$.
Let $X$ and $Y$ be the weight lattice and the root lattice of $\fg$. The order of the group $X/Y$ is det($a_{ij}$).
Let $X_+$ be the set of dominant weights. 
For any $\ga, \gb \in X$, we say $\alpha > \beta$ if $\alpha - \beta \in \BZ_+ \Pi_\fh$. This defines
a partial ordering on $X$. Let $\Phi_+ = \{\ga\in \Phi\st \ga>0\}$ be the set of positive roots.

One can define a symmetric bilinear form ($\cdot | \cdot$) on $\fh^*$ in 
the following way. Multiply the $i$-th row of ($a_{ij}$) by $d_i \in \{1, 2, 3\}$ 
such that ($d_i a_{ij}$) is a symmetric matrix with $d_i=1,\il$ or $|\{d_i\}|=2$.
Set ($\alpha_i | \alpha_j$) $= d_i a_{ij}$. 
This bilinear form is non-degenerate and is proportional to the dual of the Killing form restricted on $\fh$. 
Set $d = \max_{1\leq i \leq \ell}\{d_i\}$. Let $\ga_0$ be the highest {\em short} root with respect to $\Pi_\fh$.
The Coxeter number $h = (\rho | \ga_0)+1$ where $\rho = \half\sum_{\ga>0} \ga$. Then we have
$h = \ell +1, 2\ell, 2\ell, 2\ell-2, 12, 18, 30, 12, 6$ when
$\fg = A_\ell, B_\ell, C_\ell, D_\ell, E_6, E_7, E_8, F_4, G_2$ respectively.

The fundamental dominant weights $\lambda_i,\ \il$ are defined by 
$2(\lambda_i | \alpha_j)/(\alpha_j | \alpha_j) = \dij$. Here $\dij$ is the Kronecker symbol. 
Then $X$ and $X_+$ are $\BZ$-span and $\BZ_+$-span of $\lambda_i's$.

Let $\Vg$ be the category of finite-dimensional
$\fg$-modules over $\BC$. Then 
$\Vg$ is a semisimple tensor category whose simple objects are parameterized by the dominant weights, i.e. for each
dominant weight $\gl$ there exists a unique simple $\fg$-module $\gVl$ of highest weight $\gl$ and all simple
objects in $\Vg$ are of this form. 

The Weyl group $W$ acts on $X$ and $Y$ naturally. For each integer $m > 0$, the affine 
Weyl group $W_m$ is a group of linear transformations of $\fh^*$ generated by the Weyl group $W$ 
and the reflection $s_m$ along the hyperplane $\digamma\!_m = \{x\in \fh^*\st (x | \ga_0) = m\}$.
There is another $W_m$ action on $\fh^*$, called the dot action, defined
by $w . x = w(x + \rho) - \rho$ where the action on the right-hand side is the natural one. 
Set 
$$
C_m = \{x\in \fh^*\st (x + \rho | \ga_0) < m, 0<(x+\rho | \ga_i) \textrm{ for } \il \}.
$$
A fundamental domain of the dot action of $W_m$ on $\fh^*$ is $\bar C_m$, the topological closure of $C_m$. 
It is well known that $W_m = W\ltimes m Y$.

%%%%%%%%%%%%%%%%%%%%%%%%%%%%%%%%%%%%%%%%%%%%%%%%
\subsection{Quantum groups}\label{sec:quantum}
Notation in this subsection is consistent with that of \cite{j} except that we substitute $q$ there by $v$.
Recall that the quantum integer
$$
[n]_i = \frac{v_i^n - v_i^{-n}}{v_i - v_i\oo}
$$
for $n\in\BZ$ where $v_i = v^{d_i}$, the quantum factorial
$[n]^!_i = \prod_{k=1}^{n} [k]_i$ for $n\ge 0$, and the quantum binomial coefficient
$$
\left[ {a}\atop{b} \right]_i = [a]_i[a-1]_i\cdots [a-b+1]_i/[b]^!_i
$$
with $a, b \in \BZ$ and $b\ge 0$ (with the convention that $\left[ {a}\atop{0} \right]_i = 1)$.

The quantum group $\Uv = \Uv(\fg)$ associated to $\fg$ is a non-commutative $\qv$-algebra generated by
$E_i, F_i, K_i$ and $K_i\oo$ (for $\il$) with the following relations.
For $i, j = 1, \ldots, \ell$,
$$\ba{ll}
K_i K_i\oo = 1 = K_i\oo K_i, & K_i K_j = K_j K_i, \\
K_i E_j K_i\oo = v^{(\ga_i | \ga_j)}E_j, & K_i F_j K_i\oo = v^{-(\ga_i | \ga_j)} F_j,
\ea$$ 
$$
E_i F_j -F_j E_i = \dij\frac{K_i - K_i\oo}{v_i - v_i\oo}
$$
and if $i \neq j$
$$
\sum_{k=0}^{1-a_{ij}} (-1)^k \left[ {1-a_{ij}}\atop{k}\right]_i
E_i^{1-a_{ij}-k} E_j E_i^k = 0,
$$
$$
\sum_{k=0}^{1-a_{ij}} (-1)^k \left[ {1-a_{ij}}\atop{k}\right]_i
F_i^{1-a_{ij}-k} F_j F_i^k = 0.
$$
For $\ga = \sum_i k_i \ga_i \in Y$, we denote $\prod_i K_i^{k_i}$ by $K_\ga$.

A $\Uv$-module $M$ is said to have type 1 if it admits a weight space decomposition
\be\label{eq:weight}
M = \bigoplus_{\gl\in X} M^\gl
\ee
where $M^\gl = \{x\in M\st K_\ga(x) = v^{(\ga | \gl)} x\}$. If $M^\gl\neq 0$ then $\gl$ is called a weight of $M$.
A type 1 $\Uv$-module is said to have highest weight $\gl$ 
if all of its weights are less than $\gl$, with respect
to the partial ordering on $X$ defined in \secref{sec:Lie}.

One can put a Hopf algebra structure on $\Uv$ with coproduct $\Delta_{\Uv}$, antipode
$S_{\Uv}$, and counit $\epsilon_{\Uv}$ defined by
$$\ba{lll}
\Delta_{\Uv}(E_i) = E_i\otimes 1 + K_i\otimes E_i, & \epsilon_{\Uv}(E_i) = 0, & S_{\Uv}(E_i) = -K_i\oo E_i,\\
\Delta_{\Uv}(F_i) = F_i\otimes K_i\oo + 1\otimes F_i, & \epsilon_{\Uv}(F_i) = 0, & S_{\Uv}(F_i) = -F_i K_i,\\
\Delta_{\Uv}(K_i) = K_i\otimes K_i, & \epsilon_{\Uv}(K_i) = 1, & S_{\Uv}(K_i) = K_i\oo.
\ea$$
Therefore the category $\Vv$ of finite-dimensional $\Uv$-modules of type 1 is a tensor category with duality. 
The tensor product and the duality are the usual ones. The unit object $\1 = \1_{\Vv}$ is the one dimensional representation $_v\!V_0$. The coevaluation and the evaluation are
$$\ba{ll}
b_V : \1 \to V\otimes V^* \quad & \textrm{defined by} \quad b_V(1) = \sum_i e_i \otimes e^i, \\
d_V : V^*\otimes V \to \1 \quad & \textrm{defined by} \quad d_V(f\otimes a) = f(a),
\ea$$
where $\{e_i\}$ and $\{e^i\}$ are bases of $V$ and $V^*$ with $e^i(e_j) = \dij$. Let $\vVl$, $\gl\in X_+$,
be the module of highest weight $\gl$ in $\Vv$ with $\dim M^\gl = 1$. Recall that $w_0$ is the longest element of $W$.
We have
\begin{prop}[cf. \cite{lu} proposition 6.3.6]\label{easy_red}
$\Vv$ is semisimple and $\vVl$, $\lambda \in X_+$ deplete all simple objects in $\Vv$. Furthermore $\vVl^*\cong {\vV_{-w_0(\gl)}}$.
\end{prop}

Recall that $\cA = \zv$. Let $\UA$ be the $\cA$-subalgebra of $\Uv$ generated by elements
$E_i^{(n)}, F_i^{(n)}, K_i, K_i\oo$ for $\il$ and $n\in\BZ$, $n>0$. The Hopf algebra structure on $\Uv$ can be 
restricted to $\UA$. Let $\bU_\cA^+, \bU_\cA^-$ and $\bU_\cA^0$ be the subalgebras of $\UA$ generated by
$\{E_i^{(n)}\}, \{F_i^{(n)}\}$ and $\{K_i, K_i\oo, \left[ {K_i; c}\atop{t}\right]\st c, t \in \BZ,
t >0 \}$ respectively. Here
$$
\left[ {K_i; c}\atop{t}\right] = \prod_{s = 1}^t 
{\frac {K_i v_i^{c+1-s} - K_i\oo v_i^{-c-1+s}}{v_i^s - v_i^{-s}}}.
$$
Then $\UA = \bU_\cA^- \bU_\cA^0 \bU_\cA^+$.

The subalgebra $\bU_\cA^-$ is a free $\cA$-module. 
It has two bases, the PBW basis {\bf P} and the canonical basis {\bf B}. 
The canonical basis is unique while the PBW basis is unique up to the choice of 
a reduced expression of $w_0$.
The definitions of these bases are quite involved and we refer the 
readers to \cite{j} for details.

%%%%%%%%%%%%%%%%%%%%%%%%%%%%%%%%%%%%%%%%%%%%%
\subsection{TQFTs based on a ribbon category}\label{tqfts}
This subsection is a recap of section 2.2 of \cite{chenle2}. See also IV of \cite{tu1} and 
section 3.3 of \cite{bk}. Manifolds are always orientable and smooth. Maps between manifolds are always smooth. Non-zero (tangent or normal) vectors are equivalent up to scalar multiple. We fix a ribbon category $\cR$ in this
subsection.

%%%%%%%%%%%%%%%%%%%%%%%%%%%%%%%%%%%%%%%%%%%%
\subsubsection{Ribbon graphs}\label{rg}
A band is a homeomorphic image of $[0,1]$ with a non-zero normal field. The images of $0$ and $1$ are called the bases of the band. An annulus is a homeomorphic image of $S^1$ with a non-zero normal field. A band or an annulus is oriented if it is equipped with a non-zero tangent field. A coupon is a homeomorphic image of $[0,1]\times[0,1]$. The images of $[0,1]\times 0$ and $[0,1]\times 1$ are called the bottom and the top of the coupon respectively. Coupons are always oriented.

Let $M$ be a 3-manifold. A ribbon graph $\Omega$ in $M$ is a union of oriented
bands, annuli, and coupons embedded in $M$ such that the bases of bands are on $\partial M$ or the 
bottoms/tops of coupons.

Let $\Omega$ be a ribbon graph. An $\cR$-coloring of $\Omega$ is an assignment to each band and annulus of 
$\Omega$ an arbitrary object of $\cR$, and to each coupon of $\Omega$ a morphism of $\cR$.
A ribbon graph $\Omega$ together with a $\cR$-coloring $\gm$ is called a $\cR$-colored 
ribbon graph, denoted $\Omega(\gm)$. A ribbon graph is partially $\cR$-colored if some 
bands and/or annuli and/or coupons are $\cR$-colored.

%%%%%%%%%%%%%%%%%%%%%%%%%%%%%%%%%%%%%%%%%%
\subsubsection{Extended surfaces}\label{es}
An $\cR$-mark on a closed surface $\Gamma$ is a point $p$ on $\Gamma$ associated
with a triple $(t, V, \nu)$, where $t$ (direction of the mark) is a non-zero 
tangent vector at $p$, $V$ (label of the mark) is an arbitrary object from $\cR$ 
and $\nu$ (sign of the mark) is $+$ or $-$. An $\cR$-extended surface (or $e$-surface for short)
is a closed oriented surface $\Gamma$ together with a finite set of $\cR$-marks on it and a decomposable
Lagrangian subspace of $H_1(\Gamma, \BQ)$. If $\Gamma$ is an $e$-surface we denote by $-\Gamma$ another $e$-surface
obtained from $\Gamma$ by reversing the orientation of $\Gamma$, keeping the Lagrangian subspace
and, for every $\cR$-mark on $\Gamma$, changing its sign while keeping its label and 
direction unchanged. The empty surface is considered an $e$-surface with no $\cR$-marks.

An $e$-homeomorphism between two $e$-surfaces is a homeomorphism between the underlying surfaces
respecting the extended structure.

%%%%%%%%%%%%%%%%%%%%%%%%%%%%%%%%%%%%%%%%%
\subsubsection{Extended cobordisms}\label{e3m}
An $\cR$-extended 3-manifold is a triple $(M, \gW, w)$ that consists of an oriented 3-manifold $M$, an integer $w$
(weight of $M$) and an $\cR$-colored ribbon graph $\gW$ sitting in it. The boundary of $M$ is an $e$-surface 
that is compatible with $\gW$.

An $\cR$-extended cobordism (or $e$-cobordism for short) is a triple $(M, \Gamma, \Lambda)$ where $M$ is an
$\cR$-extended 3-manifold and $\partial M = (-\Gamma)\sqcup \Lambda$ (as an $e$-surface). An $e$-homeomorphism 
between two $e$-cobordisms is a homeomorphism between the underlying 3-manifolds respecting the extended 
structures on 3-manifolds and their boundaries.

If $(M, \Gamma, \Lambda)$ and $(M', \Gamma', \Lambda')$ are two $e$-cobordisms and there is an $e$-homeomorphism $f : \Lambda \to \Gamma'$. One can glue these two $e$-cobordisms along $f$ to get a new $e$-cobordism 
$(M\cup_f M', \Gamma, \Lambda')$. Here $M\cup_f M'$ is an $\cR$-extended 3-manifold with an 
$\cR$-colored ribbon graph (obtained by gluing ribbon graphs in $M$ and $M'$) sitting inside 
and weight computed as in IV.9.1 of \cite{tu1}.

%%%%%%%%%%%%%%%%%%%%%%%%%%%%%%%%%%%%%%%%%%%%%%%%%%
\subsubsection{TQFTs based on $\cR$}\label{rtqft}
Let $K$ be a commutative ring with unit and $\Mod(K)$ be the category of projective $K$-modules. 
Let $\CR$ be the category of $e$-surfaces and $e$-homeomorphisms. A topological quantum 
field theory (TQFT) based on $\cR$ with ground ring $K$ is a pair $(\cT, \tau) = (\cT_\cR, \tau_\cR)$.
Here $\cT : \CR \to \Mod(K)$ is a modular functor, cf. III.1.2 of \cite{tu1}, based on $\cR$ with ground ring $K$ 
and $\tau$ assigns to every $e$-cobordism $(M, \Gamma, \Lambda)$ a $K$-homomorphism
$$
\tau(M) = \tau(M, \Gamma, \Lambda) : \cT(\Gamma) \to \cT(\Lambda)
$$
that satisfies the naturality, multiplicativity, functorality and normalization axioms.

Let $(\cT, \tau)$ be a TQFT based on $\cR$ with ground ring $K$. For an $e$-cobordism $(M, \emptyset, \Gamma)$ 
one has $\tau(M): K \to \cT(\Gamma)$. Denote $\tau(M)(1)\in \cT(\Gamma)$ by $[M]$, called a vacuum state in
$\cT(\Gamma)$. The TQFT is called non-degenerate if for every $e$-surface $\Gamma$, the module $\cT(\Gamma)$ is
generated over $K$ by vacuum states.

An $e$-cobordism with closed underlying 3-manifold is called closed. The image of a closed 
$e$-cobordism $(M, \emptyset, \emptyset)$ under a TQFT is a linear map from $K$ to $K$. 
Therefore $[M] \in K$. The following definition is due to Gilmer \cite{g2}.

\bdfn\label{dfn:Almost}
Let $D$ be a Dedekind domain contained in $K$ and $\cE$ be an element in $K$. A TQFT $(\cT, \tau)$ is called almost ($D$-)integral if $\cE[M]$ is in $D$ for every closed connected $e$-cobordism $(M, \emptyset, \emptyset)$.
\edfn

%%%%%%%%%%%%%%%%%%%%%%%%%%%%%%%%%%%%%%%%%%%%%%
\subsection{Twisted modules}\label{twisted}
Let $K$ be a commutative ring and $A$ be a $K$-algebra. Let $V$ be an $A$-module and $f : A\to A$ be a $K$-algebra endomorphism. One can define another $A$-module structure on $V$ induced by $f$. We denote $V$ with this induced $A$-module structure by $^fV$. Let $\cdot$ and $\cdot_f$ denote these two actions of $A$ on $V$ then
$$
x\cdot_f u = f(x)\cdot u, \qquad \forall\ x\in A,\ u \in V.
$$
Let $g : A \to A$ be a $K$-algebra antiendomorphism. Then one can define an $A$-module structure on 
$V^* = \Hom_K(V, K)$ through $g$. We denote $V^*$ with this induced $A$-module structure by $^gV^*$.
The action of $A$ on it is denoted by $*_g$. (Dot is used for actions on $V$ and asterisk is used for 
actions on $V^*$.) One has
$$
(x*_g\alpha)(u) = \alpha(g(x)\cdot u), \qquad \forall x\in A,\ \alpha \in V^*,\ u \in V.
$$
Proof of the following lemma is left for the readers. 

\begin{lemma}\label{lemma_ind}
Suppose $V, W$ are $A$-modules and $f : A \to A$ is a $K$-algebra automorphism.
If $g : A \to A$ is a $K$-algebra endomorphism (resp. antiendomorphism) and $\phi : {^fV} \to {^gW}$ (resp. $\phi : {^fV} \to {^gW^*}$) is an $A$-module homomorphism then $\phi$ can be considered as an $A$-module homomorphism $V\to {^{gf^{-1}}W}$ (resp. $V\to {^{gf^{-1}}W^*}$).
\end{lemma}

We adopt the following convention: if $A$ is a Hopf algebra with antipode $S_A$ and $V$ is an $A$-module 
then we will denote $^{S_A}V^*$ simply by $V^*$.

%%%%%%%%%%%%%%%%%%%%%%%%%%%%%%%%%%%%%%%%%%%%%%%%%%%%%%%
%%%%%%%%%%%%%%%%%%%%%%%%%%%%%%%%%%%%%%%%%%%%%%%%%%%%%%%
\section{Representations of quantum groups for generic $v$}\label{rep}

%%%%%%%%%%%%%%%%%%%%%%%%%%%%%%%%%%%%%%%%%%%%%%%%
\subsection{Finite-ranked $\UA$-modules}\label{generic}
%Let $\VA$ be the category of finite-ranked $\UA$-modules of type 1 (cf. equation \eqref{eq:weight}) over the ground
%ring $\cA$. Hence objects in $\VA$ are in particulr free $\cA$-modules of finite rank.
%Let's recall the following well known proposition about $\Vv$ first before we study $\VA$.
%Recall that $w_0$ is the longest element in $W$. 
%
Let $\bP$ and $\bB$ be a PBW basis and the canonical basis of $\bU^-_\cA$ respectively (cf. \secref{sec:quantum}).
Recall that $\vVl$, $\gl\in X_+$, is the simple object in $\Vv$ of highest weight $\gl$. 

Fix an element $v_\lambda \in {\vVl^\lambda}$. The $\UA$-submodule ${\AVl}$ of $\vVl$ 
generated by $v_\gl$ is free as an $\cA$-module. The canonical basis (up to scalar) $\bB_\lambda$ of $\AVl$, and
$\vVl$ of course, is the set $\{b\,v_\lambda\ne 0 \, |\, b\in \bB\}$. A PBW basis $\bP_\gl$ of ${\AVl}$
is $\{b\,v_\lambda\ne 0 \, |\, b\in \bP\}$, which is not unique up to scalar. 
Let 
$$
\AVl^* = \Hom_\cA(\AVl, \cA),
$$
which is also a $\UA$-module because the antipode $S_{\Uv}$ restricts to $\UA$. 
The following remark demonstrates that \propref{easy_red} is no longer true for $\AVl$. 

\brmk
The quantum group $\Uv(\fs\fl_2)$ is an algebra over $\qv$ generated by $E, F, K, K^{-1}$ with the relations
$$
K K\oo = K\oo K = 1,\qquad E F - F E = \frac{K - K\oo}{v - v^{-1}},
$$
$$
K E = v^2 E K,\qquad K F = v^{-2} F K.
$$
It is a Hopf algebra with antipode $S_{\Uv(\fs\fl_2)}: \Uv(\fs\fl_2) \to \Uv(\fs\fl_2)$ defined by
$$
S_{\Uv(\fs\fl_2)}(E) = -K\oo E,\quad S_{\Uv(\fs\fl_2)}(F) = -FK, \quad S_{\Uv(\fs\fl_2)}(K) = K\oo.
$$

\begin{figure}[ht]
\centering
\psfrag{1}{$e_n$}
\psfrag{2}{$e_{n-1}$}
\psfrag{3}{$e_2$}
\psfrag{4}{$e_1$}
\psfrag{5}{$e_0$}
\psfrag{6}{$[n]$}
\psfrag{7}{$[2]$}
\psfrag{8}{$[1]$}
\psfrag{9}{$[1]$}
\psfrag{A}{$[n-1]$}
\psfrag{B}{$[n]$}
\includegraphics[width=4in]{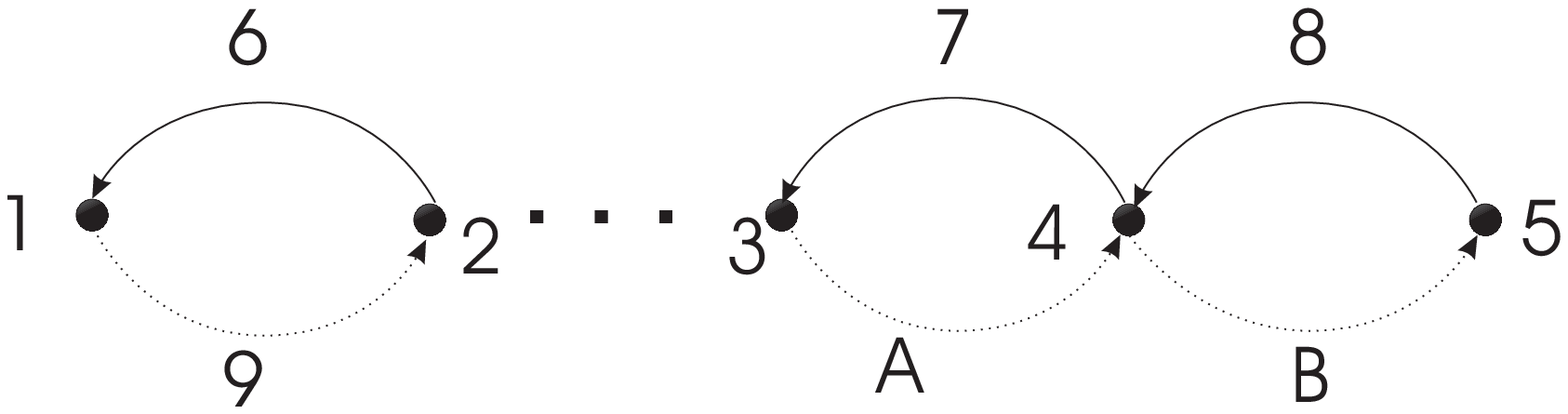}
\caption{}\label{fig:Vn}
\end{figure}

\begin{figure}[ht]
\centering
\psfrag{1}{$e_n^*$}
\psfrag{2}{$e_{n-1}^*$}
\psfrag{3}{$e_2^*$}
\psfrag{4}{$e_1^*$}
\psfrag{5}{$e_0^*$}
\psfrag{6}{$a_n[n]$}
\psfrag{7}{$a_2[2]$}
\psfrag{8}{$a_1[1]$}
\psfrag{9}{$a_n\oo[1]$}
\psfrag{A}{$a_2\oo[n-1]$}
\psfrag{B}{$a_1\oo[n]$}
\includegraphics[width=4in]{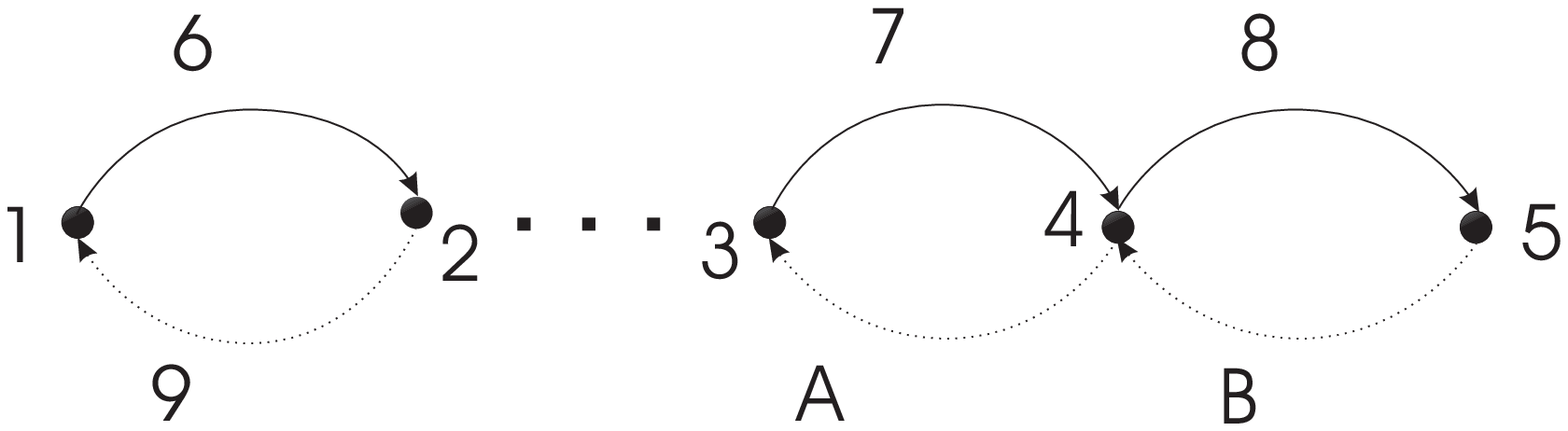}
\caption{}\label{fig:Vn*}
\end{figure}

For any integer $n \ge 0$ there exists a unique $(n+1)$-dimensional $\Uv(\fs\fl_2)$-module $V_n$
over $\qv$ with basis $\bb_n = \{e_i\}_{i = 0, \ldots, n}$ such that the actions of 
$E$ and $F$ on $V_n$ are depicted by the diagram in
\figref{fig:Vn}, 
where the action of $F$ is indicated by solid lines and the action of $E$ by dashed lines. 
Let $\bb_n^*=\{e_i^*\}_{i=0, \ldots, n}\subset V_n^*$ such that $e_i^*(e_j) = \delta_{i j}$. 
Using the antipode $S_{\Uv(\fs\fl_2)}$, one has the diagram in \figref{fig:Vn*} 
depicting the actions of $E$ and $F$ on ${V_n}^*$ where $a_k = -v^{n-2 k +2}$.
From these two diagram it is easily seen that $V_n\cong V_n^*$. It turns out that $\bb_n$ is the canonical basis of $V_n$. Hence $\bb_n$ is a basis for $\AVn$ and $\bb_n^*$ is a basis for $\AVn^*$. These two diagrams also show that $\AVn$  and $\AVn^*$ are not isomorphic as $\UA(\fs\fl_2)$-modules. Proposition~\ref{gendu} below says that in order to have isomorphism one only needs to invert quantum integers. 

\ermk

%%%%%%%%%%%%%%%%%%%%%%%%%%%%%%%%%%%%%%%%%%%%%%%%%%%
\subsection{Inverting quantum integers}\label{sec:invert}
The $\UA$-module ${\AVl}, \lambda \in X_+$, carries a unique 
bilinear form $H_\lambda$ (cf. \cite{ck}) such that:
\begin{enumerate}
\item[]
$H_\lambda(v_\lambda, v_\lambda) = 1$,
\item[]
$H_\gl(x,a y) = a H_\gl(x,y)$ for $a\in \cA$ and $x, y \in {\AVl}$,
\item[]
$H_\lambda(g\cdot x, y) = H_\lambda(x, \varphi(g)\cdot y)$ for $g \in \UA$ and $x, y \in {\AVl}$.
\end{enumerate}
Here $\varphi$ is the unique $\BZ$-algebra antiautomorphism of $\UA$ such that:
\be\label{varphi}
\varphi(E^{(n)}_i) = F^{(n)}_i,\quad \varphi(F^{(n)}_i) = E^{(n)}_i,\quad \varphi(K_i) = K_i\oo,
\quad \varphi(v) = v^{-1}.
\ee
$H_\lambda$ is often referred to as the {\em quantum Shapovalov's form}.
We have $H_\lambda({\AVl^\mu}, {\AVl^\nu}) = 0$ if $\mu \neq \nu$ and $\Ker (H_\lambda) = 0$. Let $\det(\lambda)$ denote the determinant of the matrix of $H_\lambda$ in the basis $\bP_\lambda$.
 
\begin{prop}[\cite{ck} proposition~1.9]\label{qsh}
For any $\lambda \in X_+$, $\det(\lambda)$ is a product of quantum integers.
\end{prop}

Let $\An = \cA[1/[n]^!_i]_{\il}$ and $\AnVl = {\AVl}\otimes_\cA \An$. For $\gl\in X_+$, let
$$
\AnVl^* = \Hom_{\cA_n}({\AnVl}, \cA_n).
$$
Since $\AVl$ and $\AnVl$ are free modules of same rank $\AnVl^* ={\AVl^*}{\otimes_\cA}\An$. We have

\begin{prop}\label{gendu}
For $\gl\in X_+$, the $\UA$-modules $\AnVwl$ and $\AnVl^*$ are isomorphic for large enough $n$ where $\lambda' = -w_0(\lambda)$.
\end{prop}
\begin{proof}
The bilinear form $H_\lambda: {\AVl} \otimes_\BZ {\AVl} \to \cA$ induces a 
$\UA$-linear map over ground ring $\BZ$ (cf. \secref{twisted})
\be\label{h}
H_\lambda: {\AnVl} \to {^\varphi\!\!\!\!\AnVl^*}.
\ee
Let $\bar\ : \UA \to \UA$ be the unique 
$\BZ$-algebra automorphism such that
$$
\overline{E^{(s)}_i} = E^{(s)}_i,\quad \overline{F^{(s)}_i} = F^{(n)}_i,\quad \overline{K_i} = K_i\oo,
\quad \bar v = v^{-1}.
$$
Then $\bar \varphi = \varphi \circ \bar{ }\ $ is an antiautomorphism for $\UA$ (cf. equation~(\ref{varphi}))
over $\cA$. By \lemmaref{lemma_ind}, equation~(\ref{h}) can be considered as a $\UA$-module homomorphism
$$
H_\lambda: {\bar\ \!\!\!\!\AnVl} \to {^{\bar\varphi}\!\!\!\AnVl^*}
$$
over the ground ring $\BZ$. By \propref{qsh}, $H_\lambda$ is actually an isomorphism.

There is a unique $\UA$-module isomorphism $\bar\ : {\AnVl} \to \bar\ \!\!\!\! \AnVl$, over ground ring $\BZ$, that is identity on the canonical basis $\bB_\lambda$ of $\AnVl$ (cf. \cite{lu} page 171 where $\vVl$ and $\bB_\gl$ are denoted by $\gL_\gl$ and $\bB(\gL_\gl)$ respectively). Then 
$$
\bar H_\lambda = H_\lambda \circ \bar{ } : {\AnVl} \to {^{\bar \varphi}\!\!\!\AnVl^*}
$$ 
is a $\UA$-module isomorphism over $\BZ$. Easy computation shows that $\bar H_\gl$ is an honest $\UA$-module isomorphism, i.e. over ground ring $\cA$.
There is a unique algebra automorphism $\omega$ of $\UA$ such that
$$
\omega(E^{(s)}_i) = F^{(s)}_i, \quad \omega(F^{(s)}_i) = E^{(s)}_i, \quad \omega(K_i) = K_i\oo.
$$
Lusztig proved that there exists a unique $\UA$-module isomorphism 
$$
\chi : {\AnVl} \to {^{\omega}\!\!\!\!\AnVwl}
$$
that maps $\bB_\lambda$ to $\bB_{\lambda'}$ (cf. \cite{lu} page 177). Hence
\be\label{new}
\bar H_\lambda \circ \chi^{-1} : {^{\omega}\!\!\!\!\AnVwl} \to {^{\bar \varphi}\!\!\!\!\AnVl^*}
\ee
is a $\UA$-module isomorphism. By \lemmaref{sm} below $^{S_{\Uv}\oo \bar \varphi}\!\!\!\AnVwl$ 
is isomorphic to ${^{\omega}\!\!\!\!\AnVwl}$, which in turn is isomorphic to ${^{\bar \varphi}\!\!\!\!\AnVl^*}$ 
by equation~(\ref{new}). Hence $^{S_{\Uv}\oo \bar \varphi}\!\!\!\AnVwl$ and ${^{\bar \varphi}\!\!\!\AnVl^*}$ 
are isomorphic. This isomorphism can be considered as an isomorphism between $\AnVwl$ 
and $^{S_{\Uv}}\!\!\!\AnVl^*$.
\end{proof}

\begin{lemma}\label{sm}
For any $\mu \in X_+$, ${^{S_{\Uv}\oo \bar \varphi}\!\!\!\AnVm} \cong {^{\omega}\!\!\!\!\AnVm}$.
\end{lemma}

\begin{proof}
By \lemmaref{lemma_ind}, it's enough to prove $V_1 = {^{S_{\Uv}\oo \bar \varphi \omega}\!\!\!
\AnVm} \cong V_2 = {\AnVm}$ (because $\omega^2 = \Id$).
Let $z: V_1 \to V_2$ and $z': V_2 \to V_1$ be $\UA$-homomorphisms such that $z(v_{\mu}) = v_{\mu}$ and $z'(v_{\mu}) = v_{\mu}$. They are well defined because $S_{\Uv}\oo \bar \varphi \omega$ is a $\UA$-automorphism sending $K_i$ to $K_i$. Since $_v\!V_\mu = V_1\otimes\BQ(v) = V_2\otimes\BQ(v)$ is the unique type 1 simple $\Uv$-module of highest weight $\mu$, $z\otimes \BQ(v) \circ z'\otimes \BQ(v)$ and $z'\otimes \BQ(v) \circ z\otimes \BQ(v)$ are both identity maps. Therefore $z$ and $z'$ are inverse to each other.
\end{proof}

Set $\Ainf = \cup_{n = 1}^\infty \An$ and $\AinfVl = {\AVl}\otimes_\cA \Ainf$. Because the coproduct 
$\gD_{\Uv}$ of $\Uv$ restricts to $\UA$, tensor product of $\UA$-modules is also a $\UA$-module. 
Using \propref{gendu} and the canonical basis we have

\begin{lemma}\label{highpart}
$\AinfV_{\lambda + \mu}$ is a direct summand of $\AinfV_{\lambda}\otimes {\AinfV_\mu}$ for $\gl, \gm \in X_+$.
\end{lemma}

\begin{proof}
Let $M_1 = {\AinfV_{\lambda}}$, $M_2 = {\AinfV_\mu}$, $M_3 =  M_1\otimes M_2$, 
$M_4 = {\AinfV_{\lambda + \mu}}$,
$M_5 = {\AinfV_{-w_0(\gl+\mu)}}$, and $M_6 = M_5 \otimes M_3$. Also let $M_0 = {\AinfV_0}$ be the rank 1 
$\UA$-module. Following \cite{lu}, one can define the canonical basis $\bB_i$ for $M_i$. 
Note that $\bB_3 = \bB_1 \diamondsuit \bB_2$ (resp. $\bB_6 = \bB_5 \diamondsuit \bB_3$) is {\em not}
the tensor product basis $\bB_1 \otimes \bB_2$ (resp. $\bB_5 \otimes \bB_3$).

There is a $\UA$-module monomorphism 
$$
\phi : M_4 \to M_3
$$
that carries $\bB_4$ into $\bB_3$ (cf. \cite{lu} proposition~27.1.7). 
To prove the lemma we need to find a $\UA$-module
epimorphism $\phi' : M_3 \to M_4$ such that $\phi'\circ\phi = \Id_{M_4}$.

It is easy to see that $M_6\otimes \BQ(v)$ has a unique quotient isomorphic to $M_0\otimes \BQ(v)$ considered as objects in $\Vv$. (Use weight argument and \propref{easy_red}.)
By proposition~27.2.6 in \cite{lu} one has a $\UA$-epimorphism 
$$
\pi : M_6 \to M_0
$$
carrying $\bB_6$ onto $\bB_0 \cup \{0\}$. There is a ($\UA$-linear) right coevaluation $\vartheta : \Ainf \to M_5^*\otimes M_5$ such that
\be\label{coe}
\vartheta(1) = \sum_{i=1}^m b_i^* \otimes K_{-2 \rho}(b_i).
\ee
Here $\bB_5 = \{b_i\}$ and $b_i^*(b_j) = \delta_{ij}.$ Suppose $u_4$ is of highest weight in $\bB_4$ and $b_m$ is of lowest weight in $\bB_5$. By \propref{gendu} there is a $\UA$-isomorphism $\psi : M_5^* \to M_4$ such that 
\be\label{c}
\psi(b_m^*) = c\,u_4
\ee
for some invertible $c\in \Ainf$.

Let 
$$
\phi'' = \psi\circ (\Id_{M_5^*}\otimes\pi)\circ (\vartheta\otimes\Id_{M_3}) : M_3 \to M_4.
$$
Since $M_4\otimes \BQ(v)$ is simple in the category 
$\Vv$, $(\phi''\circ\phi) \otimes \Id_{\BQ(v)} = a \Id_{M_4\otimes\BQ(v)}$ for some $a\in \BQ(v)$. 
But $\phi''$ and $\phi$ are defined over $\Ainf$ we know that $a \in \Ainf$. 

We prove that $a$ is invertible in $\Ainf$ and this will complete the proof of the lemma.
We only have to look at the image of $u_4$ under $\phi''\circ\phi$.
Let $u_3$ be the unique element of weight $\gl+\gm$ in $\bB_3$. Because $\phi$ takes canonical basis to canonical basis one has $\phi(u_4) = u_3$ and
$$
(\vartheta\otimes\Id_{M_3})\circ \phi (u_4) = \sum_{i = 1}^m b_i^* \otimes K_{-2 \rho}(b_i) \otimes u_3
$$
(cf. equation~(\ref{coe})). By theorem~24.3.3.(b) of \cite{lu} $b_m \diamondsuit u_3 = b_m \otimes u_3\in \bB_6$. 
By proposition~27.3.8 of \cite{lu} $\pi(b_m \otimes u_3) = u_0$ where $u_0$ is the unique element in $\bB_0$. 
Since $b_m$ has weight $-\gl - \gm$ one has $\pi (K_{-2 \rho}(b_m) \otimes u_3) = v^{(2 \rho | \gl + \gm)}u_0$. 
For $i\neq m$, $\pi (K_{-2 \rho}(b_i)\otimes u_3) = 0$ because $b_i\otimes u_3$ is of nonzero weight. 
Therefore $(\Id_{M_5^*}\otimes\pi)\circ (\vartheta\otimes\Id_{M_3})\circ \phi(u_4) = v^{(2 \rho | \gl + \gm)} b_m^*.$
By equation~(\ref{c}), we have $a = v^{(2 \rho | \gl + \gm)} c$, which is invertible.
\end{proof}

\begin{figure}
\centering
\begin{minipage}[c]{0.5\textwidth}
\centering
\psfrag{4}{$\vartheta$}
\includegraphics[width=1in]{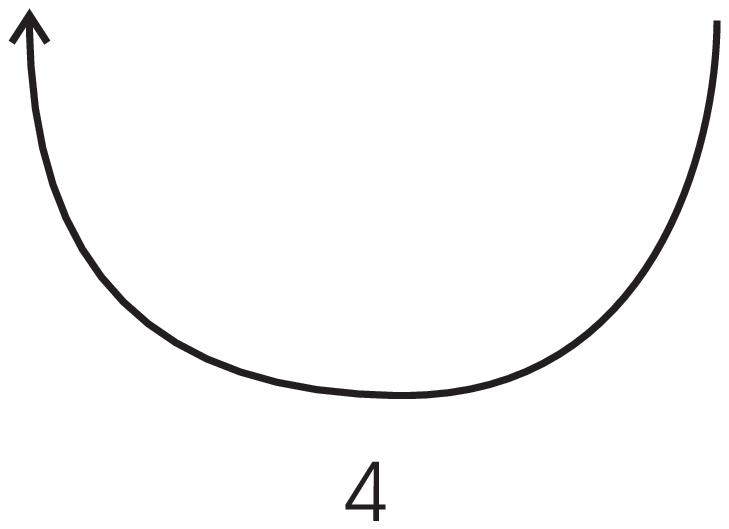}
\end{minipage}%
\begin{minipage}[c]{0.5\textwidth}
\centering
\psfrag{1}{$\psi$}
\psfrag{2}{$\pi$}
\psfrag{3}{$\phi$}
\psfrag{4}{$\AinfV_{\gl + \gm}$}
\psfrag{5}{$\AinfV_\gl$}
\psfrag{6}{$\AinfV_\gm$}
\psfrag{7}{$\AinfV_{-w_0(\gl+\gm)}$}
\psfrag{8}{$\phi''$}
\includegraphics[width=2in]{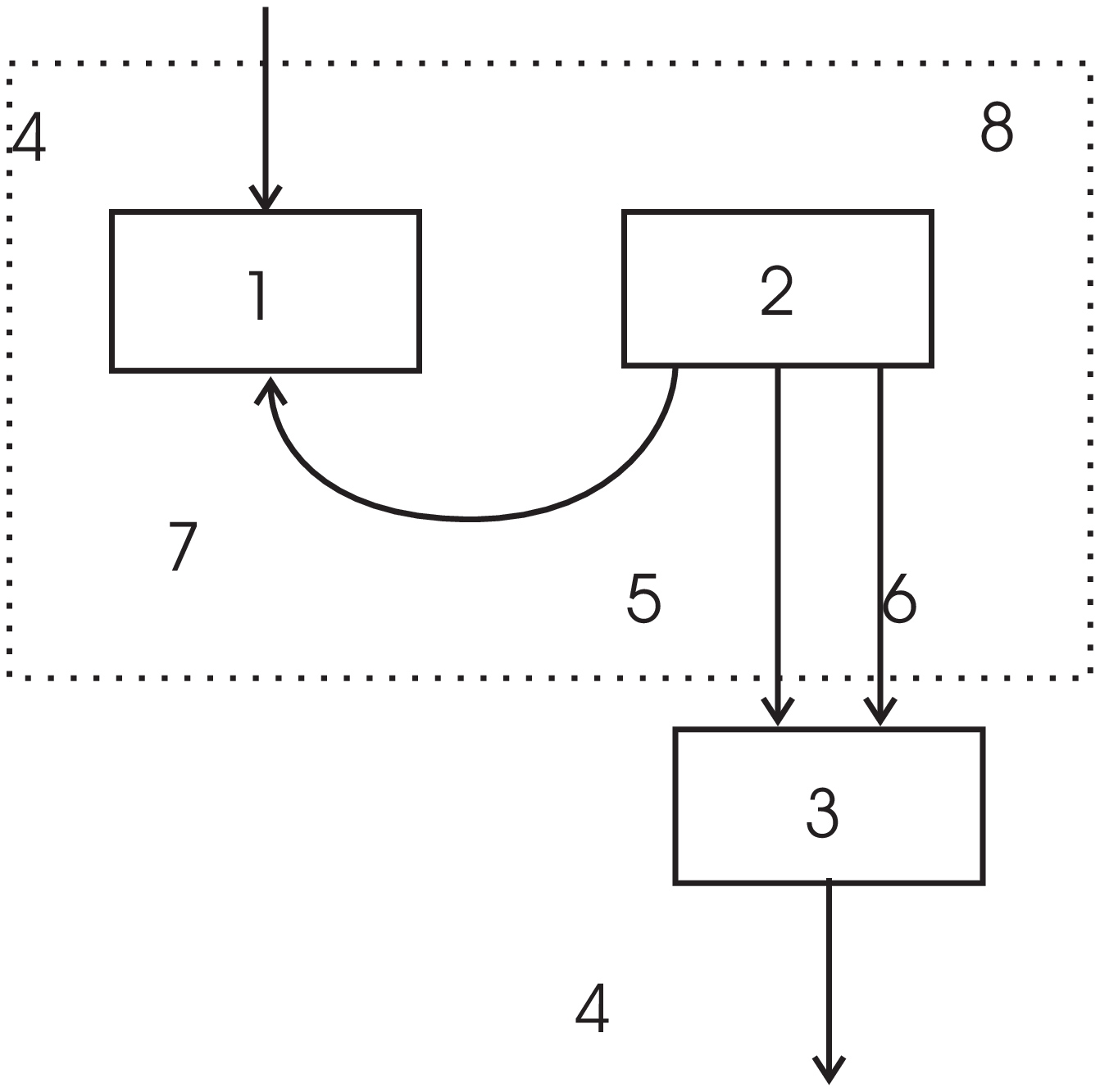}
\end{minipage}
\caption{Morphisms used in the proof of \lemmaref{highpart}}\label{fig:MorphismUsedInLemma}
\end{figure}

The finite-dimensional $\UA$-modules of type 1 form a ribbon category whose morphisms 
have ribbon graph presentations (cf. \cite{tu1}). Figure~\ref{fig:MorphismUsedInLemma}
presents the morphisms used in the proof. 

Furthermore we have

\begin{prop}\label{genten}
For $\lambda, \mu \in X_+$, there exist $\gm_1,\, \ldots,\, \gm_m \in X_+$ such that
$$
{\AinfVl}\otimes{\AinfV_\mu} = \bigoplus_{i = 1}^m {\AinfV_{\gm_i}}.
$$
\end{prop}
\begin{proof}
As we saw above, $\AinfV_{\lambda + \mu}$ is a direct summand of $\AinfV_{\lambda}\otimes {\AinfV_\mu}$. 
Let $M$ be the direct summand of $\AinfV_{\lambda}\otimes {\AinfV_\mu}$ that is complementary 
to $\AinfV_{\lambda + \mu}$. We must show that $M$ contains a direct summand isomorphic to 
$\AinfV_\nu$ for some $\nu \in X_+$. Note that $M$ is a based module (cf. \cite{lu}). 
Let $\nu$ be a highest weight of $M$. By proposition~27.1.7 in \cite{lu} $M$ contains a 
submodule isomorphic to $\AinfV_\nu$. This submodule is actually a direct summand.
This can be proved similarly as the previous lemma. The only difference is that $M$ may 
contain more than one copy of $\AinfV_\nu$ but it is easy to get around. One only has to 
be a little more careful when constructing the projection map. We leave the details for the readers.
\end{proof}

Let $\cV_\infty$ be a category of $\UA$-modules and $\UA$-morphisms over ground ring $\Ainf$.
Objects in $\cV_\infty$ are direct summands of $\AinfV_{\gm_1}\star \cdots \star {\AinfV_{\gm_m}}$
where $\star$ is either $\otimes$ or $\oplus$ and $\gm_i \in X_+$. 
For any $\gl \in X_+$, $\AinfVl$ is simple because 
even its extension to $\qv$ is simple in $\Vv$. An easy consequence 
of the previous proposition is the following.

\begin{cor}
The tensor category $\cV_\infty$ is semisimple and $\AinfVl$'s, $\gl \in X_+$, deplete all simple objects.
Furthermore ${\AinfV_\gl}^*\cong{\AinfV_{\gl'}}$ where $\gl' = -w_0(\gl)$.
\end{cor}

%%%%%%%%%%%%%%%%%%%%%%%%%%%%%%%%%%%%%%%%%%%%%%%%%
%%%%%%%%%%%%%%%%%%%%%%%%%%%%%%%%%%%%%%%%%%%%%%%%%
\section{Representations at roots of 1}\label{root1}
From now on $r$ is assumed to be an odd prime greater than the Coxeter number and
$\xi$ is a root of unity of order $r$.

%%%%%%%%%%%%%%%%%%%%%%%%%%%%%%%%%%%%%%%%%%%%%%%%
\subsection{$\Ux$-modules}\label{sec:Ux}
Let $\Ux = \UA \otimes_\cA \BZ[\xi]$ where $\BZ[\xi]$ is considered as an $\cA$-algebra 
by sending $v$ to $\xi$. Similarly we define $\Ux^-, \Ux^+$ and $\Ux^0$. 
One has $\Ux = \Ux^- \Ux^0 \Ux^+$. We study the representation theory for $\Ux$ in this section 
using results from \secref{sec:invert} and \cite{andersen}. Let $\xVl = {\AVl} 
\otimes_\cA \BZ[\xi]$ and $\xVl^* = \Hom_{\BZ[\xi]}(\xVl, \BZ[\xi])$ for $\gl \in X_+$. 
They are $\Ux$-modules and are free of finite rank as $\BZ[\xi]$-modules.
They will be the building blocks of our construction. 
The next two lemmas concern their duality and tensor product.

\begin{lemma}\label{1du}
For $\gl \in \bar C_r \cap X_+$ and $\lambda' = -w_0(\lambda)$ the $\Ux$-modules $\xVl$ and $\xVwl^*$ are isomorphic.
\end{lemma}

Note that $\gl' \in \bar C_r \cap X_+$ if and only if
$\gl \in \bar C_r \cap X_+$ where $\bar C_r$ is defined in \secref{sec:Lie}. 
Let $\CVl = {\xVl} \otimes \BC$, $\gl \in X_+$. It is well known that $\CVl$ 
is a simple $\UC$-module for $\gl \in \bar C_r \cap X_+$. Note that the image 
of the quantum integer $[n]_i$ in $\zx$ is invertible 
if $r$ is not a factor of $n$ and the image is 0 otherwise.

\begin{proof}
The quantum Shapovalov's form $H_\lambda$ on $\AVl$ (cf. \secref{generic}) induces a bilinear form on $\xVl$.
Reading the proof of \propref{gendu} one realizes that it's enough to show that $\det(\gl)$ 
contains no factor of $[sr]$ for some integer $s$. Assume there was one. Then $\det(\gl) = 0$. 
Since $H_\lambda$ can also induce a bilinear form on $\CVl$ this would imply that $\CVl$ has 
a nontrivial submodule, namely the kernel of $H_\lambda$. We get a contradiction.
\end{proof}

\begin{lemma}\label{1tensor}
If $\gl, \mu, \gl + \mu \in X_+\cap \bar C_r$ then $\xVl\otimes {\xVm} = \oplus_{i = 1}^s {_\xi\! V_{\mu_i}}$ with
$\mu_i \in X_+ \cap \bar C_r$.
\end{lemma}
\begin{proof}
The proof goes like those of \lemmaref{highpart} and \propref{genten}. Remember we had to prove that the element $a = v^{(2 \rho | \gl + \gm)} c$ is invertible, which follows from the fact that $c$ is invertible. Situation here is similar. We have to prove one element is invertible, which can be proved by using the proof of \lemmaref{1du} and the discussion before it.
\end{proof}

%%%%%%%%%%%%%%%%%%%%%%%%%%%%%%%%%%%%%%%%%%%%%%%%%%%%%%%%%%%%%%%%%%%%%%%%
\subsection{Negligible objects and morphisms}\label{sec:negligible}
Let $\Vx'$ be a full subcategory of $\Vx$, the category of finite-ranked $\Ux$-modules of type 1.
Every object in $\Vx'$ is a direct summand of $\xVmo\star \cdots \star \xVmm$ 
where $\star$ is either $\otimes$ or $\oplus$ and $\gm_i \in Y_+ \cap C_r$ where $Y_+ = Y\cap X_+$. 
Morphisms in $\cV'_\xi$ are $\Ux$-morphisms. (We use the root lattice $Y$ instead of $X$ because
(1) it's required in the case $\fg = B_\ell$ and (2) it gives modular category.)

An object $V$ in $\cV'_\xi$ is called negligible if the quantum trace 
$\tr_q(f) = \tr(K_{2\rho}f) = 0$ for all $f\in\End_{\Vx'}(V)$.
A morphism $f:V\to W$ in $\Vx'$ is called negligible if it factors through a negligible object.
(Compare \cite{kirillov} section~3.)

We record some basic facts about negligible objects and morphisms.

\begin{lemma}\label{lm:facts}
In $\Vx'$ we have
\begin{enumerate}
	\item $M\otimes N$ is negligible if $M$ is,
	\item $\xVl$ is negligible for every $\gl \in \digamma\!_r\cap Y_+$ (cf. \secref{sec:Lie}),
	\item $\xVl$, $\gl \in C_r\cap Y_+$, is not a direct summand of a negligible module,
	\item direct summand of a negligible object is negligible,
	\item direct sum of negligible objects is negligible,
	\item $f$, $f^*$, $f g$, $g f$, $f\otimes g$ and $g\otimes f$ are negligible if $f$ is,
	\item $f+g$ is negligible if $f$ and $g$ are.
\end{enumerate}
\end{lemma}

\begin{proof}
One can prove (1) and (6) using graphical calculus. We prove (2). Any $f\in\End_{\Vx'}(\xVl)$ induces
$f\otimes 1\in \End_{\VC}(\CVl)$. We see that $\tr_q(f) = 0$ because $\tr_q(f\otimes 1) = 0$
(cf. \cite{kirillov}). We prove (3). If $\xVl$, $\gl \in C_r\cap Y_+$, was a direct summand of a
negligible module $M$ then $\xVl\otimes\BC$ would be a direct summand 
of $M\otimes\BC$. But Andersen proved this can not be true (cf. \cite{kirillov} again). Assertion (4) is obvious.
Assertion (5) follows from $\tr_q(f) =\tr_q(\pi_1\circ f\circ \iota_1)+\tr_q(\pi_2\circ f\circ \iota_2)$
for $f\in\End_{\Vx'}(V_1\oplus V_2)$ where $\pi_i$ and $\iota_i$, $i = 1, 2$, are standard projection and injection.
Assertion (7) follows from the universal property of the direct sum and (5).
\end{proof}

Let
\bea\label{eq:bound}
m(A_\ell) = \ell + 1, & m(B_\ell) = 2\ell, & m(C_\ell) = 3\ell -1,\\
m(D_\ell) = 3\ell - 6, & m(E_6) = 14, & m(E_7) = 21. \nonumber
\eea
Note that $m(\fg)\ge h$. The following proposition is the key to prove \propref{prop:tensor}.

\begin{prop}\label{prop:decomposition}
If $\fg$ is not of type $E_8, F_4$ or $G_2$ and $r$ is a prime bigger than $m(\fg)$ then
for any object $V$ in $\Vx'$ we have
$$
V = \bigoplus_{i}{\xVmi}\oplus Z
$$
where $\gm_i\in Y_+\cap C_r$ and $Z$ is negligible.
\end{prop}

Under the assumption $r$ is, in particular, prime to the non-zero entries of the Cartan 
matrix for $\fg$ and is bigger than the Coxeter number. This proposition will be proved in \secref{sec:proof}.

%%%%%%%%%%%%%%%%%%%%%%%%%%%%%%%%%%%%%%%%%%%%%%%%%
\subsection{Reduction to semisimple tensor categories}\label{sec:reduction}
Let $\bVx'$ be the category of $\Vx'$ quotient by negligible morphisms. Objects in $\bVx'$ are the same
as those in $\Vx'$ and
$$
\Hom_{\bVx'}(V, W) = \Hom_{\Vx'}(V, W)/\textrm{negligible morphisms}.
$$
Negligible objects in $\Vx'$ become 0 in $\bVx'$. By assertions (6) and (7) of \lemmaref{lm:facts},
$\bVx'$ is a tensor category. According to \propref{prop:decomposition}, $\bVx'$ is semisimple and is dominated
by $\xVl$, $\gl\in Y_+\cap C_r$. This proves \propref{prop:tensor}.

%%%%%%%%%%%%%%%%%%%%%%%%%%%%%%%%%%%%%%%%%%%%%%%%
%%%%%%%%%%%%%%%%%%%%%%%%%%%%%%%%%%%%%%%%%%%%%%%%
\section{Proof of \propref{prop:decomposition}}\label{sec:proof}
\begin{figure}[ht]
\centering
\psfrag{1}{$1$}
\psfrag{2}{$2$}
\psfrag{3}{$3$}
\psfrag{4}{$4$}
\psfrag{5}{$5$}
\psfrag{6}{$6$}
\psfrag{7}{$7$}
\psfrag{L}{$\ell$}
\psfrag{X}{$\ell - 1$}
\psfrag{Y}{$\ell - 2$}
\psfrag{Z}{$\ell - 3$}
\psfrag{A}{$A_\ell\quad (\ell\ge 1)$}
\psfrag{B}{$B_\ell\quad (\ell\ge 2)$}
\psfrag{C}{$C_\ell\quad (\ell\ge 3)$}
\psfrag{D}{$D_\ell\quad (\ell\ge 4)$}
\psfrag{E}{$E_6$}
\psfrag{F}{$E_7$}
\includegraphics[width=4.5in]{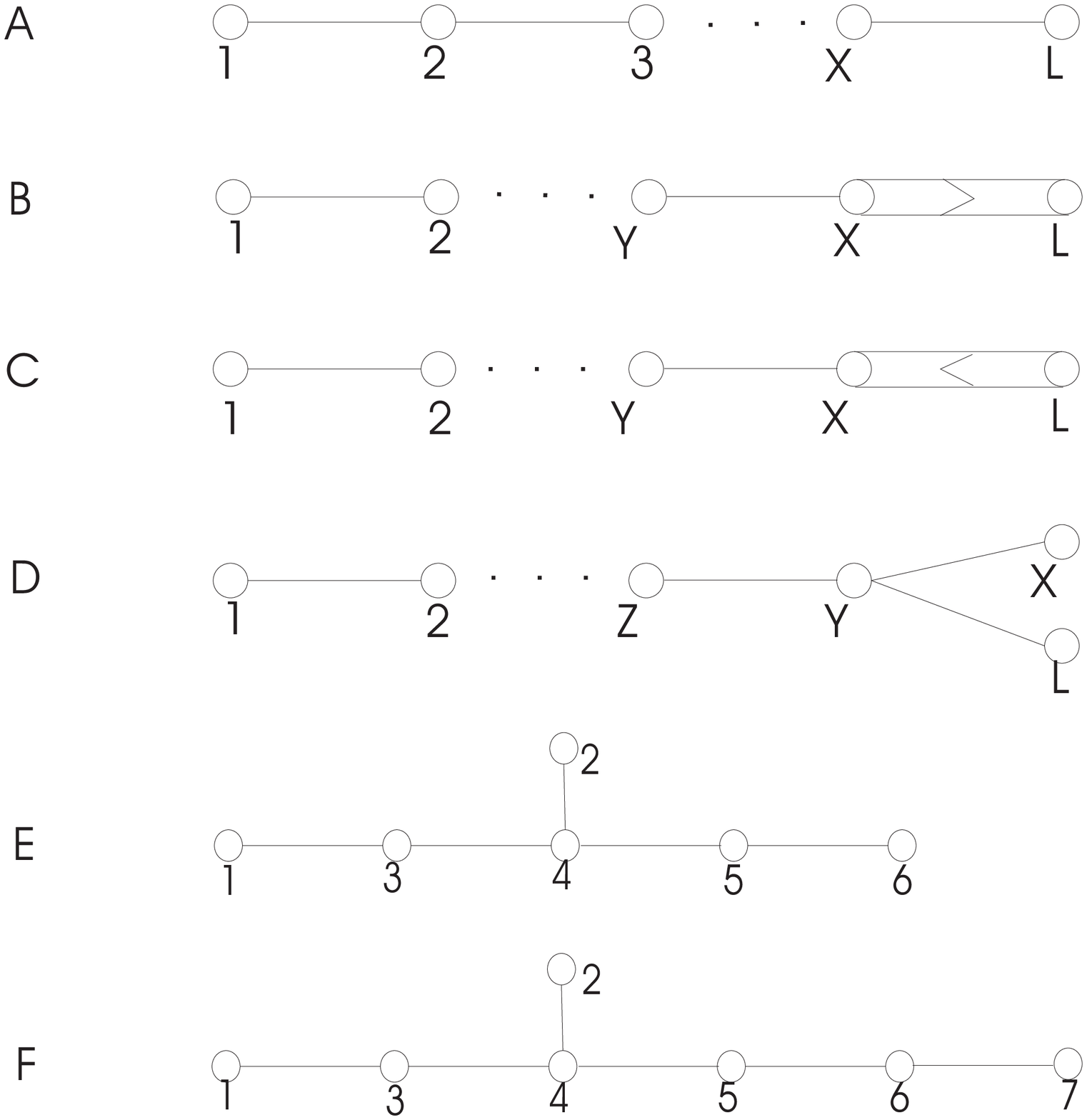}
\caption{}\label{fig:dynkin}
\end{figure}

We will make use of some techniques from \cite{gk}. The proof is done by checking different cases. In
all the cases except $B_\ell$ we need to enlarge $\Vx'$ a little bit. Let $\Vx''$ be a subcategory of $\Vx$. 
Every object in $\Vx'$ is a direct summand of $\xVmo\star \cdots \star \xVmm$ 
where $\star$ is either $\otimes$ or $\oplus$ and $\gm_i \in X_+ \cap C_r$. 
Morphisms in $\cV'_\xi$ are $\Ux$-morphisms. Obviously $\Vx'$ is a full subcategory of $\Vx''$. We can
define negligible objects and morphisms in $\Vx''$ similarly as in $\Vx'$. If we prove \propref{prop:decomposition}
in $\Vx''$ (i.e. changing $Y$ to $X$) then it is true in $\Vx'$ also. 
The Dynkin diagrams in \figref{fig:dynkin} help us identify fundamental representations.

%%%%%%%%%%%%%%%%%%%%%%%%%%%%%%%%%%%%%%%%%%%%%%%%%%
\subsection{Type $A_\ell$}\label{sec:A}
Recall that $\gl_i$, $\il$, are the fundamental weights. In this case, the highest (short) root 
$\ga_0 = \ga_1+\cdots +\ga_\ell$ (cf. \cite{h1}). We have $(\gl_i+\rho|\ga_0) = 1+\ell$. 
Hence for any $r > m(\fg) = \ell+1$, $\xVli$ is an object in $\Vx''$. 

Let $\gl = \sum_i a_i \gl_i$ be a dominant weight in $C_r\cap X_+$. By \lemmaref{1tensor}, $\xVl$ is a direct
summand of $\otimes_i({\xVli\!\!\! ^{\otimes a_i}})$. Since $\xVm$, $\gm\in X_+\cap C_r$,
is simple and direct summand of negligible object is negligible (\lemmaref{lm:facts}.(4))
it is enough to prove \propref{prop:decomposition} when $V = {\xVl\otimes\xVli}$ with $\gl\in X_+\cap C_r$.
In this case we have $(\gl+\gl_i+\rho|\ga_0) = (\gl+\rho|\ga_0) + 1 \le r$. The proposition now follows from
\lemmaref{1tensor} and \lemmaref{lm:facts}.(2).

%%%%%%%%%%%%%%%%%%%%%%%%%%%%%%%%%%%%%%%%%%%%%%%
\subsection{Type $C_\ell$}\label{sec:C}
It's known that $\xVli = {\xVlo\!\!\!\!^{\wedge i}}$, $\il$,
which is a direct summand of ${\xVlo\!\!\!\!^{\otimes i}}$ if the highest weight $i\gl_1$ is in $\bar C_r\cap X_+$
(cf. theorem 10.5.7 of \cite{wan} and \lemmaref{1tensor}). 
We have $\ga_0 = \ga_1+2(\ga_2+\cdots+\ga_{\ell-1})+\ga_\ell$ and
$(\gl_1|\ga_0) = 1$. So $(i\gl_1+\rho|\ga_0) \le (\ell\gl_1+\rho|\ga_0)= 3\ell -1 = m(\fg) < r$ 
and $i\gl_1$ is in $\bar C_r\cap X_+, \il$.
Therefore we only have to prove the proposition when $V={\xVm\otimes{\xVlo}}$ for $\gm \in C_r\cap X_+$. This can
be proved as in the case $A_\ell$.

%%%%%%%%%%%%%%%%%%%%%%%%%%%%%%%%%%%%%%%%%%%%%%%
\subsection{Type $D_\ell$}\label{sec:D}
We have $\ga_0 = \ga_1+ 2(\ga_2+\cdots+\ga_{\ell-2})+\ga_{\ell-1}+\ga_\ell$. 
For $i = 1,\ell-1$ and $\ell$ and $n = 1,2,\ldots,\ell-3$, we have 
$$
(n \gl_i +\rho |\ga_0) = n+2\ell-3\le 3\ell-6 = m(\fg)<r.
$$
According to theorem 10.5.7 of \cite{wan} and \lemmaref{1tensor}, all fundamental representations are 
direct summands of $\xVli\!\!\! ^{\otimes j}, i=1,\ell-1,\ell$ and 
$j=1,2,\ldots,\ell-3$. We can prove the proposition just like 
$A_\ell$ because $(\gl_i |\ga_0) = 1$, $i = 1, \ell-1, \ell$.

%%%%%%%%%%%%%%%%%%%%%%%%%%%%%%%%%%%%%%%%%%%%%%%
\subsection{Type $E_6$}\label{sec:E6}
We have $\ga_0 = \ga_1+2\ga_2+2\ga_3+3\ga_4+2\ga_5+\ga_6$ and $(n\gl_i+\rho|\ga_0) \le 14 = m(\fg)<r$, 
$i = 1,6$ and $n = 1,2,3$. So $\xVli, i = 1,6$ belong to $C_r\cap X_+$.
Using \cite{stembridge} and \propref{1tensor}, we see that all fundamental
representations $\xVli, 1\le i \le 6$ are direct summands of 
$\xVlo\!\!\!^{\otimes j}$ or $\xVls\!\!\!^{\otimes j}$, $j = 1,2,3$.
By \lemmaref{lm:facts}.(2) we only have to prove the proposition 
when $V={\xVm}\otimes{\xVli}, i = 1, 6$. This can
be proved as $A_\ell$ since $(\gl_i|\ga_0) = 1$ for $i = 1,6$.

%%%%%%%%%%%%%%%%%%%%%%%%%%%%%%%%%%%%%%%%%%%%%%%
\subsection{Type $E_7$}\label{sec:E7}
We have $\ga_0 = 2\ga_1+2\ga_2+3\ga_3+4\ga_4+3\ga_5+2\ga_6+\ga_7$
and $(n\gl_7+\rho|\ga_0) \le 21 = m(\fg)<r$, $n=1,2,3,4$. So
$\xVlsv$ belongs to $C_r\cap X_+$. Using \cite{stembridge} and \propref{1tensor}, we see that all fundamental
representations $\xVli, 1\le i \le 7$ are direct summands of $\xVlsv^{\otimes j}$, $j = 1,2,3,4$.
By \lemmaref{lm:facts} we only have to prove the proposition when $V={\xVm}\otimes{\xVlsv}$. This can
be proved as $A_\ell$ since $(\gl_7|\ga_0) = 1$.

%%%%%%%%%%%%%%%%%%%%%%%%%%%%%%%%%%%%%%%%%%%%%%%
\subsection{Type $B_\ell$}\label{sec:B}
We work in $\Vx'$. Let $\bb_\ell = \{\gl_1,\ldots,\gl_{\ell-1}, 2\gl_\ell\}\subset Y$. Because $|X/Y| = 2$ and
$\bb_\ell \subset Y$ we see that $\bb_\ell$ is a basis of the root lattice $Y$ over $\BZ$. 
For any $\gm \in \bb_\ell$, $(\gm|\ga_0) = 2$ where $\ga_0 = \sum_i\ga_i$, cf. \cite{h1}. 
Therefore $\xVm$, $\gm\in \bb_\ell$, is an object in $\Vx'$. For any $\gl \in C_r\cap Y_+$, $\xVl$ is a direct
summand of $\xVmo\otimes\cdots\otimes{\xVmm}$ for $\gm_i \in \bb_\ell$ by \lemmaref{1tensor}. Therefore
we only have to prove the proposition when $V = {\xVl}\otimes{\xVm}$ for $\gl\in Y_+\cap C_r$ and $\gm\in\bb_\ell$.

From the Cartan matrix we see that $(\ga_i|\ga_j)$ is always an even number. 
For any $\gl \in Y$, there is an integer $a$ such that
$$\ba{ll}
(\gl +\rho | \ga_0) & = 2 a + (\rho|\ga_0) = 2 a + \sum_{i,j = 1}^\ell (\ga_i|\gl_j)\\
 & = 2 a + \sum_{i=1}^{\ell-1}\sum_{j = 1}^\ell (\ga_i|\gl_j) + \sum_{j=1}^\ell (\ga_\ell|\gl_j)\\
 & = 2 a + \sum_{i=1}^{\ell-1}\sum_{j = 1}^\ell 2\dij + \sum_{j=1}^\ell \delta_{j\ell}\\
 & = 2 a + 2(\ell-1)+1 = 2(a+\ell) - 1.
\ea$$
If $\gl\in Y_+\cap C_r$ then $(\gm+\rho|\ga_0)\le r-2$ because $r$ is odd. 
Hence $\gl+\gm\in \bar C_r\cap Y_+$ for any $\gm\in\bb_\ell$ and we 
are done by \propref{1tensor} and \lemmaref{lm:facts}.(2).

%%%%%%%%%%%%%%%%%%%%%%%%%%%%%%%%%%%%%%%%%%%
%%%%%%%%%%%%%%%%%%%%%%%%%%%%%%%%%%%%%%%%%%%
\section{Integral TQFTs}\label{tqft}
We will follow \cite{tu1} to construct a modular category from the braided tensor category $\bVx'$.
Then we will follow \cite{g2} to get integral TQFTs.
The construction of the TQFT is similar to the one in \cite{chenle2}.

Let $\fg$ be a complex simple Lie algebra not of type $E_8, F_4, G_2$ and $r$ be an odd prime greater than
$m(\fg)$. We fix a set of dominant simple objects 
$$
A = A(\bVx') = \{\xVl\st \gl \in Y_+\cap C_r\}.
$$
We identify $A$ with the set $\{\gl\in Y_+\cap C_r\}$.

%%%%%%%%%%%%%%%%%%%%%%%%%%%%%%%%%%%%%%%%%%%
\subsection{Modular categories}\label{s}
It is well known that $\bVx'$ is a ribbon category, cf. \cite{chenle2, le1}. 
In order to make it into a modular category, we only have to make the $S$-matrix 
$$
S=S(\bVx') = (s_{\gl,\gm})_{\gl,\gm\in A}
$$ 
invertible.

Let $\BC[X]$ be the algebra over $\BC$ generated by $v^\gl$, $\gl\in X$. The addition in $\BC[X]$ 
is formal and the multiplication is given by $v^\gl v^\gm = v^{\gl + \gm}$. One can map $\BC[X]$ 
onto a subalgebra of $\Hom_\BC(\Ux^0, \BC)$ through 
$$
v^\gl(K_i) = \psi_r\left(v^{(\ga_i|\gl)}\right),\quad v^\gl(\left[{K_i; c}\atop{t}\right]) 
= \psi_r\left(\left[{(\ga_i|\gl)/d_i+c}\atop{t}\right]_i\right),
$$
where $\psi_r:\cA\to \zx$ is a ring homomorphism sending $v$ to $\xi$.
Let $\delta\in \BC[X]$ be the Weyl denominator
\be\label{Weyl_deno}
\delta = \prod_{\ga\in \Phi_+}(v^{\ga/2} - v^{-\ga/2}) = \sum_{w\in W}\sn(w)v^{w(\rho)},
\ee
where $\sn(w) = (-1)^{\textrm{length\ of\ } w}$. Let
$$
\cD^2=\cD_{\bVx'}^2 = \sum_{\gl \in Y_+\cap C_r} \dim_q(\xVl)^2,
$$
where $\dim_q(\xVl) = \tr_q(\Id_{\xVl})$ is the quantum dimension of $\xVl$.
Recall that $w_0$ is the longest element in $W$
and $\ell$ is the dimension of $\fh$.

\begin{lemma}\label{lm:S}
Let $r>m(\fg)$ be a prime and $\fg$ be a simple Lie algebra not of type $E_8, F_4, G_2$. Then
$\cD^2 = \frac{\sn(w_0)r^\ell}{\gd(K_{2\rho})^2}$ and $\det(S)^2 = \pm\cD^{2|Y_+\cap C_r|}$.
\end{lemma}

Note that $\gd(K_{2\rho}) \ne 0$ because $r$ is odd and greater than $h$.
This lemma is proved in the proof of theorem~3.3 in \cite{le2} where $v^2$ is denoted $q$. 
See also theorem~3.3.20 in \cite{bk} where $v$ is denoted $q$.
Our assumption on $r$ guarantees that it is coprime with $d \det(a_{ij})$. 

If we extend the ground $\zx$ of $\bVx'$ to $\zxr$ then the $S$-matrix will be invertible. For some technical reason
we also want to include a square root $\cD$ of $\cD^2$ in the ground ring. 
Let $\zeta$ be a root of unity of order
\be\label{eq:zeta}
O(\zeta) = \left\{\ba{ll}
r & \textrm{if\ $\ell$\ is\ even,\ $r$\ is\ arbitrary\ and\ }\ \sn(w_0)=1,\\
r & \textrm{if\ $\ell$\ is\ odd\ and\ } \sn(w_0)r\equiv 1 \mod 4,\\
4r & \textrm{otherwise.}
\ea\right.
\ee
For elements $x, y \in \BZ[\xi]$ we say that $x \sim y$ if $x = z y$ for some invertible $z$ in $\BZ[\xi]$.

\begin{lemma}\label{DinK}
$\cD \in \zzr$.
\end{lemma}

\begin{proof}
By equation~(\ref{Weyl_deno}), $\gd(K_{2\rho})\sim (\xi - 1)^{|\Phi_+|}$. 
Since $r\sim (\xi-1)^{r-1}$, $\frac{1}{\gd(K_{2\rho})}$ is in $\zzr$. 
It remains to notice that $\sqrt{\sn(w_0) r^\ell} \in \BZ[\zeta]$.
\end{proof}

Let $\zVl = {\xVl}\otimes\zzr$. Let $\Vz'$ be a category of $\Ux$-modules and $\Ux$-morphisms. Objects in $\Vz'$ are
direct summands of $\zVmo\star \cdots \star \zVmm$ 
where $\star$ is either $\otimes$ or $\oplus$ and $\gm_i \in Y_+ \cap C_r$.
Let $\bVz'$ be the category $\Vz'$ quotient by negligible morphisms. Clearly
$\bVz'$ is a modular category dominated by $A=\{{\zVl}\st\gl\in Y_+\cap C_r\}$ and $\cD$ is in the ground ring.

There is a unique functor $e_r : \bVx'\to\bVz'$ respecting the ribbon structure with $e_r(\xVl)={\zVl}$
and $e_r(f) = f\otimes\Id_{\zzr}$.

%%%%%%%%%%%%%%%%%%%%%%%%%%%%%%%%%%%%%%%%%%%%%%%
\subsection{A TQFT based on $\bVx'$}\label{sec:tqft1}
We define a TQFT based on the ribbon category $\bVx'$ (cf. section 3.2 of \cite{chenle2} and \secref{tqfts}). 
In the rest of this paper $e$-surfaces and $e$-cobordisms will mean $\bVx'$-extended surfaces and 
$\bVx'$-extended cobordisms unless otherwise specified.

%%%%%%%%%%%%%%%%%%%%%%%%%%%%%%%%%%%%%%
\subsubsection{Parametrization of $e$-surfaces}\label{pe}
\begin{figure}
\centering
\psfrag{1}{$(V_1, \nu_1)$}
\psfrag{2}{$(V_m, \nu_m)$}
\psfrag{3}{$p_1$}
\psfrag{4}{$p_m$}
\includegraphics[width=4in]{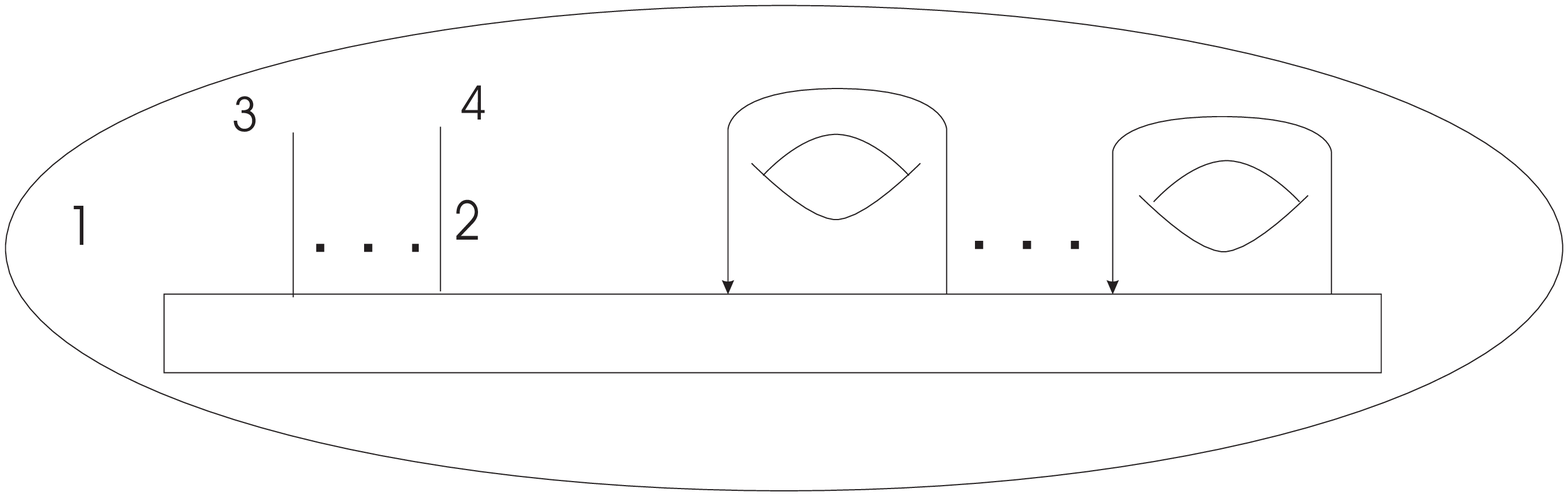}
\caption{The standard handlebody $\cH_t$}\label{hb}
\end{figure}

A $\bVx'$-extended type (or type for short) is a tuple $(g; (V_1, \nu_1), \ldots, (V_m, \nu_m))$ 
where $g\ge -1$ is an integer, $V_i$'s are arbitrary objects from $\bVx'$ and $\nu_i$'s belong to $\{+, -\}$. 
If $g = -1$ then $m = 0$. We construct a standard handlebody $\cH_t$ for every type 
$t = (g; (V_1, \nu_1), \ldots, (V_m, \nu_m))$ shown in \figref{hb}.
It is a genus $g$ handlebody standardly embedded in $\BR^3$ with a partially $\bVx'$-colored ribbon 
graph $R_t$ sitting in it (cf. \secref{rg}). Here we assume the empty space is a handlebody of genus $-1$. 
The ribbon graph $R_t$ consists of a coupon (the narrow rectangle near the bottom), $m$ vertical bands 
(with $m$ bases on the coupon and $m$ bases $p_i, \im$ on $\partial \cH_t$) and $g$ half-circled 
bands (oriented to the left with bases on the coupon). The $m$ vertical bands are oriented and colored 
according to $t$ such that the $i$-th vertical band is colored by $V_i$ and oriented up (resp. down) if $\nu_i$ is $-$
(resp. $+$). The normal vectors on the bands point toward the reader. Recall that the normal vector $n_i$ at 
$p_i$ is tangent to $\partial \cH_t$. The boundary of $\cH_t$ is an $e$-surface, called the standard $e$-surface
of type $t$ and denoted $\Sigma_t$. The $\bVx'$-marks $p_i, \im$ are associated with $(n_i, V_i, \nu_i)$.
The Lagrangian subspace is the kernel of the map $H_1(\partial \cH_t, \BQ)\to H_1(\cH_t, \BQ)$ 
induced by the inclusion. Note that $\cH_t$ is generally {\em not} an $e$-cobordism because $R_t$ 
is only partially colored.

For any connected $e$-surface $\Gamma$ let 
$$
\gp(\Gamma) = \{p\ |\ p: \Sigma_t \to \Gamma \textrm{ is an $e$-homeomorphism}\}
$$
be the set of all parametrizations of $\Gamma$ up to $e$-isotopy (in the obvious sense). Note that $\Gamma$ may be parametrized by more than one standard surfaces because the $\bVx'$-marks on $\Gamma$ are not ordered. For any $e$-surface $\Gamma$, $\gp(\Gamma)$ is not empty (IV.6.4.2 of \cite{tu1}).

%%%%%%%%%%%%%%%%%%%%%%%%%%%%%%%%%%%%%%%%%%%%%%%%%
\subsubsection{The modular functor}\label{mf}
For any integer $a$ let $A^a$ be the set $(C_r\cap Y_+)^a$ if $g>0$
and the one-element set $\{\bullet\}$ otherwise. Let $\Lambda_\bullet = \zzr$ and for $\gl\in A^g$ with $g>0$ let
$$
\Lambda_\gl = \bigotimes_{i = 1}^g (\zVli\otimes \zVli^*)
$$
where $\gl_i$ is the $i$-th coordinate of $\gl$ and the tensor products are taken over $\zzr$. Note that $\Lambda_\gl$ is an object in $\bVz'$. 

For any ($\bVx'$-extended) type $t = (g; (V_1, \nu_1), \ldots, (V_m, \nu_m))$ and $\gl \in A^g$ let
$$
T_{t, \gl} = \Hom_{\bVz'}(\zzr, (\bigotimes_{i = 1}^m e_r(V_i)^{\nu_i}) \otimes \Lambda_\gl)
$$
where the tensor products are taken over $\zzr$. Here $e_r(V_i)^+ = e_r(V_i)$ and 
$e_r(V_i)^- = e_r(V_i)^*$ for $\im$.

Let $T_t = \oplus_{\gl\in A^g} T_{t, \gl}$. Suppose $\Gamma$ is a connected $e$-surface of type $t$. 
For a parametrization $p\in\gp(\Gamma)$, set $\cT_p(\Gamma) = T_t$. 
For two parametrizations $p, p' \in \gp(\Gamma)$ define an isomorphism 
$\varphi(p, p') : \cT_p(\Gamma) \to \cT_{p'}(\Gamma)$ of 
$\zzr$-modules\footnote{The actual formulas of $\varphi(p, p')$ will not be used in the sequel. 
Interested readers can find them in IV.6.3 and IV.6.4.2 of \cite{tu1}.}. 
Identify the vector spaces $\{\cT_p(\Gamma)\}_{p\in\gp(\Gamma)}$ along the isomorphisms 
$\{\varphi(p, p')\}_{(p,p')}$. The resulting vector spaces $\cT(\Gamma)$ depends only on 
$\Gamma$. For any parametrization $p: \Sigma_t \to \Gamma$, $\cT(\Gamma)$ is canonically 
isomorphic to $\cT_p(\Gamma)$. Denote this isomorphism by $p_\sharp$.

Let $\Gamma$ and $\Gamma'$ be connected $e$-surfaces. For an $e$-homeomorphism 
$f: \Gamma\to \Gamma'$ we define $\cT(f): \cT(\Gamma)\to \cT(\Gamma')$ as follows. 
Pick any parametrization $p: \Sigma_t \to \Gamma$. Then $f p$ is a parametrization 
of $\Gamma'$. Set $\cT(f) = (f p)_\sharp (p_\sharp)^{-1}$ which does not depend on 
the parametrization we choose (cf. IV.6.3.1 of \cite{tu1}).

For non-connected $e$-surfaces we can do the above componentwisely and then form tensor product. Recall that $\cC(\bVx')$ is the category of $e$-surfaces and $e$-homeomorphisms. We have a modular functor $\cT : \cC(\bVx') \to \Mod(\zzr)$ (cf. IV.6.3.3 of \cite{tu1}).

\brmk
Unlike \cite{tu1}, which uses only a modular category in the construction, 
we start with the ribbon category $\bVx'$ in the construction of the $e$-surfaces 
and $e$-cobordisms and use the modular category $\bVz'$ later when we define the modular 
functor $\cT$. This modification ensures that the TQFT constructed in the next subsection is 
almost integral.
\ermk

{\noindent{\em Remark.\ }}
We quote several results from \cite{tu1} in this and the next subsections. They are valid in our construction because the places where we use these results are the places where we are using the modular category $\bVz'$.

%%%%%%%%%%%%%%%%%%%%%%%%%%%%%%%%%%%%%%%%%%%%%
\subsubsection{An almost integral TQFT}\label{anitqft}
Because $\bVz'$ is a modular category one can follow IV.9 in \cite{tu1} to get a TQFT $(\cT^e, \tau^w)$
based on $\bVz'$ with ground ring $\zzr$. Any $\bVx'$-extended cobordism $(M, \Gamma, \Lambda)$ induces 
a $\bVz'$-extended cobordism $(\check M, \check \Gamma, \check \Lambda)$. One just needs to replace any 
object $V \in \bVx'$ by $e_r(V)$ and any morphism $f\in \bVx'$ by $e_r(f)$. Hence one has a $\zzr$-linear map 
$$
\tau^w(\check M) : \cT^e(\check \Gamma)\to\cT^e(\check \Lambda).
$$
From the definition of the modular functor $\cT$ given in \secref{mf} it is easy to see that $\cT(Z) = \cT^e(\check Z)$ for any $\bVx'$-extended surface $Z$. Set $\tau(M) = \tau^w(\check M)$.

\begin{prop}\label{thm_main}
Suppose $\fg$ is not of type $E_8, F_4, G_2$ and $r\ge m(\fg)$ is an odd prime.
The pair $(\cT, \tau)$ is a non-degenerate TQFT based on $\bVx'$ with ground ring $\zzr$.
Furthermore, $(\cT, \tau)$ is almost $\BZ[\zeta]$-integral.
\end{prop}

This proposition can be proved as lemma 3.1 and theorem 3.2 in \cite{chenle2}. We will
sketch the proof for the almost integrality in \secref{sec:surgery}.

%%%%%%%%%%%%%%%%%%%%%%%%%%%%%%%%%%%%%%%%%
\subsection{Integral TQFTs}\label{int_tqft}
Gilmer showed in \cite{g2} that when $\fg = \fs\fl_2$ the TQFT constructed in \cite{bhmv} can be 
restricted to a sub-cobordism such that the modules of surfaces are free over some ring of algebraic 
integers. We prove that the same can be done for the TQFT $(\cT,\tau)$ constructed in \secref{sec:tqft1}.
We again fix a simple Lie algebra $\fg$ not of type $E_8, F_4, G_2$ and an odd prime $r>m(\fg)$. 

%%%%%%%%%%%%%%%%%%%%%%%%%%%%%%%%%%%%%%%
\subsubsection{Restricted cobordisms}\label{sec:restricted}
A cobordism $(M, \Sigma', \Sigma)$ is said to be targeted if the $0$-th relative Betti number
$\beta_0(M, \Sigma)$ is zero. We say $[M] = \tau(M)(1)\in \cT(M)$ is a targeted vacuum state if 
$(M,\emptyset, \Sigma)$ is targeted. A cobordism $(M, \Sigma', \Sigma)$ is said to be even if
\be\label{even}
w(M) \equiv \beta_0(\Sigma) + \beta_1(\Sigma)/2 + \beta_0(K(M)) + \beta_1(K(M)) \mod 2,
\ee
where $K(M)$ is a homologically canonical closed 3-manifold, which is equal to $M$ if $M$ 
is closed (cf. \cite{g2}). We say $[M]$ is an even vacuum state if
$(M,\emptyset, \Sigma)$ is even. 

For any $e$-surface $\Sigma$, let $\cS(\Sigma)$ (resp. $\cT_+(\Sigma)$, $\cS_+(\Sigma)$) be the submodule 
of $\cT(\Sigma)$ generated by targeted (resp. even, even targeted) vacuum states over $\zz$ (resp. $\zxr$, $\zx$).

\begin{prop}\label{prop:integral}
Given the data above, for any targeted (resp. even, even targeted) cobordism $(M, \Sigma', \Sigma)$, 
$$
\tau(M): \cS(\Sigma')\to \cS(\Sigma),\quad (resp.\ \cT_+(\Sigma') \to \cT_+(\Sigma),\ 
\cS_+(\Sigma') \to \cS_+(\Sigma)).
$$
For any $e$-surface $\Sigma$, $\cT(\Sigma)$ (resp. $\cT_+(\Sigma)$, $\cS(\Sigma)$, $\cS_+(\Sigma)$) 
is a free $\zzr$- (resp. $\zxr$-, $\zz$-, $\zx$-) module of finite rank.
\end{prop}

We will prove this proposition in \secref{proof}. It implies \propref{prop:TQFT}
and contains theorem 5 of \cite{g2} because $\sn(w_0) = -1$ and $\ell = 1$ when $\fg = \fs\fl_2$.

%%%%%%%%%%%%%%%%%%%%%%%%%%%%%%%%%%%%%%
\subsubsection{Surgery formula}\label{sec:surgery}
Let $J = J^\fg$ be the functor from $\bVx'$-colored ribbon graphs to $\bVx'$ that respects the ribbon structure,
cf. I.2.5 of \cite{tu1}. The functor $J^{\fs\fl_2}$ is usually called the colored Jones polynomial. 
Let $\gW\sqcup L$ be a ribbon graph in $S^3$ where $\gW$ is $\bVx'$-colored and $L$ is an uncolored framed link 
of $m$ components. Recall that for any integer $a\ge 0$, $A^a = (C_r\cap Y_+)^a$. 
For $\gm = (\gm_1, \ldots, \gm_m) \in A^m$, denote by $L(\gm)$ the $\bVx'$-extended framed link with the $i$-th component colored by $\xVmi$. Let $U^m$ be the trivial link of $m$ components. Define
$$
Q_{\gW\sqcup L(\mu)} = J_{U^m(\mu)} J_{\gW\sqcup L(\mu)}
$$
and
\be\label{eq:F}
F_{(L, \gW)} = \sum_{\gm\in A^m} Q_{\gW\sqcup L(\gm)}.
\ee

Let $(M, \gW, w)$ be a closed connected $\bVx'$-extended 3-manifold.
We first recall how to compute $[M]$ from a surgery description of $M$. Suppose that $M$ is the result 
of surgery on a framed link $L \subset S^3$. Then there exists a $\bVx'$-colored ribbon graph 
$\gW'$ in $S^3$ disjoint from $L$ such that $(M, \gW)$ is the result of $(S^3, \gW'\sqcup L)$ 
surgery on $L$. By slight abuse of notation we will write $\gW$ for $\gW'$. Let $U_+$ and $U_-$ 
be the trivial knots with +1 and $-1$ framing respectively. It's known that
$F_-F_+ = \cD^2$ where $F_\pm = F_{(U_\pm, \emptyset)}$ as in equation \eqref{eq:F},
cf. \cite{kirillov}. Let $\kappa$ be the square root of $F_-/F_+$ such that $\kappa \cD = F_-$.

According to IV.9.2 and II.2.2 of \cite{tu1}
\bea
[M] & = & \tau (M, \gW, w)(1) = F_-^\gs 
F_{(L, \gW)}\cD^{-(\gs + m +1)} \kappa^w \nonumber \\
& = &  F_{(L, \Omega)}
\cD^{-(m+1)} \kappa^{w+\sigma} \nonumber\\
& = & \frac{F_{(L, \Omega)}}{F_-^{\sigma_- + \beta_1}
F_+^{\sigma_+}}
F_-^{-1} \kappa^{\beta_1 + w+1}\nonumber
\eea
where $\beta_1$ is the first Betti number of $M$ and $\gs, \gs_+$ and $\gs_-$ are the signature, the number of positive eigenvalues and the number of negative eigenvalues of the linking matrix of $L$ respectively. Recall that $\sigma = \sigma_+ - \sigma_-$ and $m = \sigma_+ + \sigma_- + \beta_1$.

\begin{lemma}\label{alm}
If $(M, \gW, w)$ is a closed connected $\bVx'$-extended 3-manifold then $F_- [M]$ is in $\zz$ and it is in $\zx$ if $M$
is furthermore even. 
\end{lemma}

\begin{proof}
It is proved in \cite{chenle2} that 
$\frac{F_{(L, \Omega)}}{F_-^{\sigma_- + \beta_1} F_+^{\sigma_+}}$ is in $\BZ[\xi]$. Therefore we only
need to prove $\kappa^{\beta_1 + w+1}$ belongs to an appropriate ring.
By \lemmaref{DinK}, $\gk = F_-/\cD$ is in the fractional field of $\zz$. Since $\gk$ is always a root of 1 
(cf. \cite{bk}), $\gk$ is an element of $\zz$. By equation~(\ref{even}) $\beta_1 + w +1$ is an even
number if $M$ is even.
\end{proof}

This lemma implies that $(\cT, \tau)$ is $\zz$-almost integral 
and $(\cT_+, \tau)$ is $\zx$-almost integral.

\begin{cor}
For any $e$-surface $\Sigma$, $\cS(\Sigma)$ (resp. $\cS_+(\Sigma)$) is a finitely generated projective $\zz$-
(resp. $\zx$-) module and $\cS(\Sigma)\otimes_{\zz} \zzr = \cT (\Sigma)$ 
(resp. $\cS_+(\Sigma)\otimes_{\zx} \zxr = \cT_+(\Sigma)$).
\end{cor}

This can be proved using theorem~1 of \cite{g2} and \lemmaref{alm}.

%%%%%%%%%%%%%%%%%%%%%%%%%%%%%%%%%%%%%%%%%%%%%%%
\subsubsection{Proof of \propref{prop:integral}}\label{proof}
By the definition of $\cT$ in \secref{mf} and the fact that $\bVz'$ is dominated by free $\zzr$-modules
of finite rank, we see that $\cT(\Sigma)$ is free of finite rank.

The freeness of $\cT_+(\Sigma)$ can be proved similarly as proposition~7 in \cite{g2} using the freeness of
$\cT(\Sigma)$ and \propref{alm}. Note that $r$, $\xi$ and $\gz$ are denoted $p$, $A^2_p$ and $\alpha_p$ 
respectively in \cite{g2}. (Note that $\BZ[A_p]=\BZ[A^2_p]$ if $p$ is odd.)
Since $\cS_+(\Sigma)\otimes_{\zx} \zxr$ is isomorphic to
$\cT_+(\Sigma)$ one can show that $\cS_+(\Sigma)$ is a free $\zx$-module using lemma~1 
(and the discussion before it) of \cite{g2}. 
Hence $\cS(\Sigma) = \cS_+(\Sigma)\otimes_{\zx}\zz$ is also free.

%%%%%%%%%%%%%%%%%%%%%%%%%%%%%%%%%%%%%%%%%%%%%%%%
\subsubsection{Bases}\label{sec:bases}
Gilmer and Masbaum announced that explicit bases for the integral TQFT in \cite{g2} can be obtained
using the Kauffman bracket skein module. It will be interesting to know what these bases stand for in
the representation theory of quantum groups. Are they related to the canonical bases?


\begin{thebibliography}{BHMV}
\bibitem[A]{andersen}
Henning~Haahr Andersen.
\newblock Tensor products of quantized tilting modules.
\newblock {\em Comm. Math. Phys.}, 149(1):149--159, 1992.

\bibitem[BHMV]{bhmv}
C.~Blanchet, N.~Habegger, G.~Masbaum, and P.~Vogel.
\newblock Topological quantum field theories derived from the {K}auffman
  bracket.
\newblock {\em Topology}, 34(4):883--927, 1995.

\bibitem[BK]{bk}
Bojko Bakalov and Alexander Kirillov, Jr.
\newblock {\em Lectures on tensor categories and modular functors}, volume~21
  of {\em University Lecture Series}.
\newblock American Mathematical Society, Providence, RI, 2001.

\bibitem[CL]{chenle2}
Qi~Chen and Thang Le.
\newblock {Almost integral TQFTs from simple Lie algebras}.
\newblock {\em Submitted for publication}.

\bibitem[DCK]{ck}
Corrado De~Concini and Victor~G. Kac.
\newblock Representations of quantum groups at roots of {$1$}.
\newblock In {\em Operator algebras, unitary representations, enveloping
  algebras, and invariant theory (Paris, 1989)}, volume~92 of {\em Progr.
  Math.}, pages 471--506. Birkh\"auser Boston, Boston, MA, 1990.

\bibitem[G]{g2}
Patrick Gilmer.
\newblock Integrality for {TQFT}s.
\newblock arXiv:math.QA/0105059.

\bibitem[GK]{gk}
Sergei Gelfand and David Kazhdan.
\newblock Examples of tensor categories.
\newblock {\em Invent. Math.}, 109(3):595--617, 1992.

\bibitem[H]{h1}
James~E. Humphreys.
\newblock {\em Introduction to {L}ie algebras and representation theory}.
\newblock Springer-Verlag, New York, 1978.
\newblock Second printing, revised.

\bibitem[J]{j}
Jens~Carsten Jantzen.
\newblock {\em Lectures on quantum groups}.
\newblock American Mathematical Society, Providence, RI, 1996.

\bibitem[Ka]{kassel}
Christian Kassel.
\newblock {\em Quantum groups}, volume 155 of {\em Graduate Texts in
  Mathematics}.
\newblock Springer-Verlag, New York, 1995.

\bibitem[Ki]{kirillov}
Alexander~A. Kirillov, Jr.
\newblock On an inner product in modular tensor categories.
\newblock {\em J. Amer. Math. Soc.}, 9(4):1135--1169, 1996.

\bibitem[Le1]{le1}
Thang T.~Q. Le.
\newblock Integrality and symmetry of quantum link invariants.
\newblock {\em Duke Math. J.}, 102(2):273--306, 2000.

\bibitem[Le2]{le2}
Thang T.~Q. Le.
\newblock Quantum invariants of 3-manifolds: integrality, splitting, and
  perturbative expansion.
\newblock In {\em Proceedings of the Pacific Institute for the Mathematical
  Sciences Workshop ``Invariants of Three-Manifolds'' (Calgary, AB, 1999)},
  volume 127, no. 1-2, pages 125--152, 2003.

\bibitem[Lu]{lu}
George Lusztig.
\newblock {\em Introduction to quantum groups}.
\newblock Birkh\"auser Boston Inc., Boston, MA, 1993.

\bibitem[S]{stembridge}
John Stembridge.
\newblock {\em The coxeter and weyl Packages}.
\newblock www.math.lsa.umich.edu/~jrs/maple.html.

\bibitem[T]{tu1}
V.~G. Turaev.
\newblock {\em Quantum invariants of knots and 3-manifolds}.
\newblock Walter de Gruyter \& Co., Berlin, 1994.

\bibitem[W]{wan}
Zhe~Xian Wan.
\newblock {\em Lie algebras}.
\newblock Pergamon Press, Oxford, 1975.
\newblock Translated from the Chinese by Che Young Lee, International Series of
  Monographs in Pure and Applied Mathematics, Vol. 104.

\end{thebibliography}
\end{document}